\documentclass[american,english]{article}
\usepackage[LGR,T1]{fontenc}
\usepackage{textcomp}
\usepackage[latin9]{inputenc}
\usepackage{cprotect}
\usepackage{float}
\usepackage{amsmath}
\usepackage{amssymb}
\usepackage{graphicx}
\usepackage{geometry} % Added for better margins if needed, or rely on default
\usepackage{booktabs} % For professional tables
\usepackage{caption}
\usepackage{subcaption} % Better handling of subfigures if needed
\usepackage{authblk} 

\usepackage{algorithm}
\usepackage{algpseudocode}
\usepackage{caption} % Required to remove the algorithm number

\makeatletter

\providecommand{\tabularnewline}{\\}
\floatstyle{ruled}
\newfloat{algorithm}{tbp}{loa}
\providecommand{\algorithmname}{Algorithm}
\floatname{algorithm}{\protect\algorithmname}

\makeatother

\usepackage{babel}
\addto\captionsamerican{\renewcommand{\algorithmname}{Algorithm}}

\begin{document}
\title{Computational modeling of crack-tip fields in transversely isotropic strain-limiting solids subjected to piecewise linear slope loads}
\author[1]{S. M. Mallikarjunaiah} 
\author[2]{Saugata Ghosh}
\author[3]{Dambaru Bhatta}

\affil[1]{Department of Mathematics \& Statistics\\
Texas A\&M University-Corpus Christi\\
Corpus Christi, Texas 78412-5825, USA \\
Email: m.muddamallappa@tamucc.edu (Corresponding Author)}
\affil[2,3]{School of Mathematical \& Statistical Sciences\\
The University of Texas - Rio Grande Valley\\
Edinburg, Texas 78539, USA\\
Email:  saugata.ghosh01@utrgv.edu, dambaru.bhatta@utrgv.edu}

\date{}
\maketitle
\begin{abstract}
Crack-tip fields within a transversely isotropic strain-limiting elastic body are investigated under the influence of piecewise linear slope boundary loads. The mechanical response is characterized via a nonlinear constitutive framework relating the Cauchy stress to the linearized strain, by which non-physical strain singularities at the crack tip are eliminated. The governing system is formulated as a quasi-linear elliptic boundary value problem in terms of the displacement field and is solved utilizing a continuous Galerkin finite element method coupled with a Picard linearization scheme. Boundary conditions are prescribed such that the vertical displacement varies piecewise linearly along the top and bottom edges, exhibiting opposite slopes on each half of the boundary. Numerical results are derived for two distinct fiber orientations. It is demonstrated that piecewise slope loads provide a flexible and realistic configuration for elucidating the interplay between strain-limiting behavior and crack-tip mechanics in transversely isotropic media.\bigskip{}

\end{abstract}

\section{Introduction}
Many modern engineering structures are composed of advanced materials that possess inherent directional characteristics. Prominent examples include fiber-reinforced polymer composites, laminated plates, wood, biological tissues such as bone, and layered geological formations, all of which exhibit transversely isotropic elastic behavior \cite{kulvait2013}. In such materials, the underlying microstructure dictates that cracks often initiate and propagate along preferred directions, resulting in fracture patterns that differ significantly from those observed in isotropic solids. Consequently, the accurate characterization of stress and strain fields in the vicinity of cracks is essential for predicting failure mechanisms, assessing structural integrity, and designing reliable engineering components.

Classical linear elastic fracture mechanics (LEFM) \cite{Anderson2005,Inglis1913,love2013treatise} offers a convenient and widely adopted framework for analyzing such crack-tip behavior. However, the theory is founded on a linear constitutive relation between stress and strain, which leads to a well-known physical inconsistency: in the presence of sharp geometric features---such as crack tips, corners, and notches---the predicted strains become unbounded. The occurrence of infinite or excessively large strains is not physically realistic and may compromise the validity of LEFM at small length scales or in soft materials where the strain response saturates. These shortcomings have motivated the development of alternative constitutive models in which the strain remains bounded even under arbitrarily large stresses.

To address these limitations, strain-limiting theories inspired by Rajagopal's implicit constitutive framework \cite{rajagopal2003implicit,rajagopal2007elasticity,rajagopal2011non} have emerged as a robust alternative. In these models \cite{bulivcek2014elastic,bustamante2009some,itou2018states,kulvait2013}, the relationship between the Cauchy stress and the linearized strain is nonlinear and is specifically designed such that the strain asymptotically approaches a finite upper bound as the stress grows \cite{rajagopal2011modeling}. This concept has been successfully extended to transversely isotropic solids \cite{gou2023finite,gou2015modeling,itou2018states,Mallikarjunaiah2015,yoon2022CNSNS} and employed to study crack-tip fields under various loading conditions. These studies demonstrate that the strain singularities inherent to classical elasticity can be eliminated while retaining a smooth stress field and a mathematically well-posed boundary value problem. In the limit of vanishing strain-limiting parameters, this nonlinear model consistently reduces to the standard linear elastic response, providing a direct basis for comparison.

Despite these theoretical advances, most existing numerical studies on strain-limiting crack problems have focused on relatively simple loading scenarios, such as uniform tensile loading or smoothly varying traction conditions. In practical engineering applications, however, structural components are frequently subjected to non-uniform boundary conditions arising from partial contact, localized actuation, or complex distributed loads. Such loading configurations can induce intricate deformation patterns, including bending-like modes, asymmetric opening, or mixed-mode behavior near the crack tip. Therefore, incorporating realistic, non-uniform loads into the analysis is crucial for understanding crack-tip mechanics in situations that more closely resemble in-service conditions.

The present work investigates a transversely isotropic body containing an edge crack and subjected to a specific class of non-uniform boundary conditions: piecewise linear slope-type displacements prescribed on the top and bottom edges of a rectangular plate. Mathematically, the mechanical response is modeled within a small-strain setting, where the linearized strain tensor is related to the Cauchy stress through a special hyperelastic constitutive relationship. The combination of the equilibrium equation and this constitutive relation yields a quasi-linear elliptic boundary value problem for the displacement field. The resulting problem is solved numerically using a continuous Galerkin finite element method coupled with a Picard iterative linearization scheme.

The primary objective of this study is to quantify how the crack-tip fields---specifically stress, strain, and strain energy density---are influenced by the presence of piecewise slope loads and by the material parameters of the strain-limiting model. We explicitly contrast the predictions of the strain-limiting model with those of its linear elastic counterpart. The results demonstrate that piecewise slope boundary conditions provide a versatile loading configuration for probing the interplay between material anisotropy and strain-limiting behavior, highlighting significant differences in the structure and magnitude of the crack-tip fields compared to linear theory.

The remainder of the paper is organized as follows. Section 2 introduces the kinematic setting and the strain-limiting constitutive model for transversely isotropic bodies. In Section 3, we formulate the boundary value problem, derive its weak form, and discuss the conditions ensuring well-posedness. Numerical experiments, including parametric studies regarding the strain-limiting parameters and comparisons between distinct fiber orientations, are presented and discussed in Section 4. Finally, conclusions and potential directions for future work are summarized in Section 5.

\section{Preliminaries and Constitutive Model}

\subsection{Small strain kinematics }

We consider an open, bounded domain $D\subset\mathbb{R}^{2}$ describing
a transversely isotropic rectangular plate with a single edge crack.
The boundary is assumed to be Lipschitz continuous, is denoted by
$\partial D$ and is partitioned as 
\begin{equation}
\partial D=\Gamma_{D}\cup\Gamma_{N},\qquad\Gamma_{D}\cap\Gamma_{N}=\emptyset,\label{eq:mf1}
\end{equation}
where $\Gamma_{D}$ denotes the nonempty Dirichlet boundary, i.e.
the part of the boundary where displacements are prescribed and $\Gamma_{N}$
represents Neumann boundary, i.e., the part where traction boundary
conditions are provided. A one-dimensional subset $\Gamma_{c}\subset D$
represents a pre-existing crack that splits the domain. To describe
the motion of the body, we assume that a point $\boldsymbol{X}$ in
the reference configuration is mapped to $\boldsymbol{x}$ in the
deformed configuration. The displacement field is denoted by $\boldsymbol{u}:D\to\mathbb{R}^{2}$
and described as $\boldsymbol{u}=\boldsymbol{x}-\boldsymbol{X}.$
Let $\boldsymbol{F}:D\to\mathbb{R}^{2\times2}$ be the deformation
gradient, defined by the relation,
\begin{eqnarray}
\boldsymbol{F}= & \frac{\partial\boldsymbol{x}}{\partial\boldsymbol{X}} & =\boldsymbol{I}+\nabla\boldsymbol{u}\label{eq:mf2}
\end{eqnarray}
$\boldsymbol{I}$ represents the second-order identity tensor and
$\nabla$ denotes the gradient operator. The left and right Cauchy--Green
stretch tensors are denoted by $\boldsymbol{B}$ and $\boldsymbol{C}$
\cite{Gurtin1981,love2013treatise}, respectively, and are defined
as, 
\begin{eqnarray}
\boldsymbol{B}=\boldsymbol{F}\,\boldsymbol{F}^{\mathrm{T}}, & \qquad & \boldsymbol{C}=\boldsymbol{F}^{\mathrm{T}}\boldsymbol{F}.\label{eq:mf3}
\end{eqnarray}
Here the superscript $(\cdot)^{\mathrm{T}}$ is the \textit{transpose operator}. 
A measure of strain in the reference configuration is provided by
the Green--St.\ Venant strain tensor $\boldsymbol{E}:D\to\mathbb{R}^{2\times2}$,
given by
\begin{eqnarray}
\boldsymbol{E} & = & \frac{1}{2}\big(\boldsymbol{\boldsymbol{C}}-\boldsymbol{\boldsymbol{I}}\big).\label{eq:mf4}
\end{eqnarray}
On the other hand, the strain measure in the current configuration
is referred as the Almansi-Hamel strain tensor $\boldsymbol{e}:D\to\mathbb{R}^{2\times2}$
which is defined as,
\begin{eqnarray}
\boldsymbol{e} & = & \frac{1}{2}\big(\boldsymbol{\boldsymbol{I}}-\boldsymbol{\boldsymbol{B}}^{-1}\big).\label{eq:mf5}
\end{eqnarray}
Under the assumption of small displacement gradients, both the strain
measures $\boldsymbol{E}$ and $\boldsymbol{e}$ reduce to the linearized
strain tensor $\boldsymbol{\epsilon}$,

\begin{eqnarray}
\boldsymbol{\epsilon} & = & \frac{1}{2}\left[\nabla\boldsymbol{\boldsymbol{u}}+\big(\nabla\boldsymbol{\boldsymbol{u}}\big)^{\mathrm{T}}\right],\label{eq:mf6}
\end{eqnarray}
This is the strain measure used in the subsequent analysis.

\subsection{Rajagopal's strain limiting framework}

The classical theory of linear elasticity describes an explicit relation
between the stress and strain. Although this description is very convenient,
it tends to predict unbounded strains and singularities near the crack
tips, corners, and sharp notches. The approach of the present work
is based on Rajagopal's implicit constitutive theory that addresses
the limitations of the classical theory by prescribing an implicit
relation between the Cauchy stress $\boldsymbol{\sigma}$ and the
Cauchy--Green stretch tensor $\boldsymbol{B}$ such that,

\begin{eqnarray}
\mathcal{F}(\boldsymbol{B},\boldsymbol{\sigma}) & = & \boldsymbol{0}.\label{eq:rj1}
\end{eqnarray}
A special subclass of Rajagopal\textquoteright s implicit framework
\cite{rajagopal2003implicit,rajagopal2007elasticity} that is relevant
for the present study assumes that the Cauchy--Green tensor $\boldsymbol{B}$
can be expressed explicitly in terms of Cauchy stress, namely
\begin{eqnarray}
\boldsymbol{B} & = & \mathcal{\widetilde{F}}\left(\boldsymbol{\sigma}\right).\label{eq:rj2}
\end{eqnarray}
Under the assumption of small displacement gradients, the above relation
provides linearized strain $\boldsymbol{\epsilon}$ in terms of $\sigma$,
\begin{eqnarray}
\boldsymbol{\epsilon} & = & \mathcal{\bar{F}}(\boldsymbol{\sigma}).\label{eq:rj3}
\end{eqnarray}
Here $\mathcal{\bar{F}}:Sym\left(\mathbb{R}^{2\times2}\right)\rightarrow Sym\left(\mathbb{R}^{2\times2}\right)$
is a nonlinear response function. The model described in (\ref{eq:rj3})
can be referred as \textit{strain limiting }\cite{Mallikarjunaiah2015,MalliPhD2015}
in nature if it has the distinguished feature, $\max_{\mathbf{\mathbf{\boldsymbol{\sigma}}}\in Sym}\left\Vert \mathcal{\bar{F}}(\mathbf{\mathbf{\boldsymbol{\sigma}}})\right\Vert \leq M$
for a constant $M>0$. Motivated by previous studies in this area
\cite{lee2022finite,yoon2021quasi,yoon2022MMS}, we adopt a specific
choice of the response operator $\mathcal{\bar{F}}\left(\cdot\right)$
that yields a strain limiting behavior,
\begin{eqnarray}
\mathcal{\bar{F}}\left(\boldsymbol{\mathbf{\boldsymbol{\sigma}}}\right) & = & \frac{\mathbb{K}[\mathbf{\boldsymbol{\mathbf{\boldsymbol{\sigma}}}}]}{(1+(\beta|\mathbb{K}^{\frac{1}{2}}[\mathbf{\boldsymbol{\mathbf{\boldsymbol{\sigma}}}}]|)^{\alpha})^{1/\alpha}}.\label{eq:rj4}
\end{eqnarray}
$\alpha>0$ and $\beta\geq0$ are the model parameters. $\mathbb{K}[\mathbf{\mathbf{\boldsymbol{\cdot}}}]$
represents the \foreignlanguage{american}{fourth order compliance
tensor.} Moreover $\sup_{\mathbf{\mathbf{\boldsymbol{\boldsymbol{\sigma}}}}\in Sym}\left\Vert \mathcal{\bar{F}}\left(\boldsymbol{\mathbf{\boldsymbol{\sigma}}}\right)\right\Vert \leq\frac{1}{\beta}$
, so $\beta$ specifies the maximum admissible strain limit. Also,
as the parameters $\beta\rightarrow0$ or $\alpha\rightarrow\infty$,
the nonlinear model converges to its linear counterpart. We assume
that the nonlinear response mapping $\mathcal{\bar{F}}:\mathrm{Sym}(\mathbb{R}^{2\times2})\to\mathrm{Sym}(\mathbb{R}^{2\times2})$
satisfies the conditions listed below. These properties ensure the
well-posedness of the boundary value problem.
\begin{enumerate}
\item \textbf{Uniform boundedness.} The Response operator $\mathcal{\bar{F}}$
satisfies the following inequality,
\begin{eqnarray*}
\left\Vert \mathcal{\bar{F}}(\mathbf{A})\right\Vert \le\frac{1}{\beta}.\qquad & \text{for all }\mathbf{A}\in\mathrm{Sym}(\mathbb{R}^{2\times2}).
\end{eqnarray*}
Consequently, the strain response is uniformly bounded, even when
the stress becomes large.
\item \textbf{Strict monotonicity.} For any two distinct tensors $\mathbf{A},\mathbf{B}\in\mathrm{Sym}(\mathbb{R}^{2\times2})$,
the mapping $\mathcal{\bar{F}}$ satisfies 
\begin{equation}
\big(\mathcal{\bar{F}}(\mathbf{A}_{})-\mathcal{\bar{F}}(\mathbf{B})\big):\big(\mathbf{A}-\mathbf{B}\big)>0.
\end{equation}
Thus $\mathcal{\bar{F}}$ is strictly monotone in the sense of the
Frobenius inner product, a property that rules out multiple solutions
to the associated variational problem and underlies uniqueness.
\item \textbf{Lipschitz continuity.} There exists a constant $\tilde{k}_{1}>0$
such that 
\begin{equation}
\left\Vert \mathcal{\bar{F}}(\mathbf{A})-\mathcal{\bar{F}}(\mathbf{B})\right\Vert \le\tilde{k}_{1}\,\left\Vert \mathbf{A}-\mathbf{B}\right\Vert \qquad\text{for all }\mathbf{A},\mathbf{B}\in\mathrm{Sym}(\mathbb{R}^{2\times2}).
\end{equation}
In other words, $\mathcal{\bar{F}}$ is globally Lipschitz on the
space of symmetric tensors, so small perturbations in the argument
produce proportionally small changes in the response.
\item \textbf{Coercivity.} There exists a constant $\tilde{k}_{2}>0$ such
that 
\begin{equation}
\left\langle \mathcal{\bar{F}}(\mathbf{A}):\mathbf{A}\right\rangle \;\ge\;\tilde{k}_{2}\,\|\mathbf{A}\|^{2}\qquad\text{for all }\mathbf{A}\in\mathrm{Sym}(\mathbb{R}^{2\times2}).
\end{equation}
\\
This property ensures the existence of a solution. It is important
to note that, the constant $\tilde{k}_{2}$ is a model-dependent quantity,
influenced by the material parameters, the particular form of the
constitutive relation, and the dimension of the problem. In this sense,
the mapping $\mathcal{\bar{F}}$ is coercive.
\end{enumerate}

\subsection{Constitutive relation for Transverse Isotropy}

We consider a transversely isotropic material with a preferred fiber
direction represented by a unit vector $\boldsymbol{a}\in\mathbb{R}^{2}$.
The corresponding structural tensor is 
\begin{eqnarray}
\boldsymbol{M} & = & \boldsymbol{a}\otimes\boldsymbol{a},\label{eq:TI1}
\end{eqnarray}
here, $\otimes$ denotes the tensor (dyadic) product, i.e., $\left(\boldsymbol{a}\otimes\boldsymbol{a}\right)_{ij}=a_{i}a_{j}.$
Internal anisotropy of the body is described by this symmetric, rank
one tensor $\boldsymbol{M}$ by specifying the directional dependence
of the material response. The linear elastic response of such a medium
is described by a fourth-order elasticity tensor $\mathbb{E}$ , which
is basically the inverse of previously described fourth-order compliance
tensor $\mathbb{K}$. The elasticity tensor $\mathbb{E}$ for the
transversely isotropic material can be expressed as,
\begin{eqnarray}
\mathbb{E}[\boldsymbol{\epsilon}] & = & 2\mu\left(\boldsymbol{\epsilon}+\frac{\lambda}{2\mu}\,\mathrm{tr}(\boldsymbol{\epsilon})\,\boldsymbol{I}+\overline{\gamma}\,(\boldsymbol{\epsilon}:\boldsymbol{M})\,\boldsymbol{M}\right).\label{eq:TI2}
\end{eqnarray}
Here $\mu>0$ and $\lambda>0$ are the Lam\'e's parameters
and $\overline{\gamma}$ is an additional anisotropy parameter governing
stiffness along the preferred direction, $\boldsymbol{I}$ is the
identity tensor, and ``$:$'' denotes the inner product,
\begin{eqnarray*}
\boldsymbol{A}:\boldsymbol{B} & = & \sum_{i,j}A_{ij}B_{ij}.
\end{eqnarray*}
Notably, the constitutive relation (\ref{eq:rj3}) is invertible for
small strains $\boldsymbol{\epsilon}$ and  within strain limiting
setting, the Cauchy stress $\boldsymbol{\sigma}$ can be expressed
as a nonlinear function of the linear strain $\boldsymbol{\epsilon}$.
The relation is expressed as follows,
\begin{eqnarray}
\boldsymbol{\boldsymbol{\sigma}}\left(\boldsymbol{\boldsymbol{\epsilon}}\right) & = & \frac{\mathbb{E}[\boldsymbol{\epsilon}]}{\left\{ 1-\left(\beta\left|\mathbb{E}^{1/2}[\boldsymbol{\boldsymbol{\epsilon}}]\right|\right)^{\alpha}\right\} ^{1/\alpha}}.\label{eq:TI3}
\end{eqnarray}
The boundary value problem discussed in the current work is based
upon the hyperelastic formulation mentioned above. The primary emphasis
of the study is to assess how the crack tip fields -- stress, strain
and strain-energy density -- predicted by the strain-limiting model
differ from those obtained from the classical linear elastic model,
which can be recovered by choosing $\beta=0$ in (\ref{eq:TI3}).

\section{Boundary Value Problem and Variational Formulation}

\subsection{Strong form of the boundary value problem}

In this section, we formulate the boundary value problem for a transversely
isotropic, strain-limiting elastic body containing a pre-existing
crack and subjected to piecewise slope-type boundary conditions on
the top and bottom boundaries. The equilibrium equation combined with
the hyperelastic constitutive relation (\ref{eq:TI3}) yields a quasi
- linear elliptic partial differential equation with displacement
field $\boldsymbol{u}=\left(u_{1},u_{2}\right)$ as the unknown. The
strong form of the boundary value problem becomes,
\begin{eqnarray}
-\nabla\cdot\boldsymbol{\boldsymbol{\sigma}} & = & \boldsymbol{\mathbf{f}}\qquad\mathrm{in}\;\:D,\nonumber \\
\boldsymbol{u} & = & \boldsymbol{u}_{D}\qquad\text{on }\Gamma_{D},\nonumber \\
\boldsymbol{\sigma n} & = & \boldsymbol{\boldsymbol{g}}\qquad\mathrm{on}\;\;\Gamma_{N},\label{eq:bvp1}
\end{eqnarray}
$\boldsymbol{u}_{D}$ is the prescribed displacement on the Dirichlet
boundary $\Gamma_{D}$, and $\boldsymbol{g}$ is the prescribed traction
on the Neumann boundary $\Gamma_{N}$. The body force per unit volume
is denoted by $\boldsymbol{\mathbf{f}}\in\left(L^{2}\left(D\right)\right)^{2}.$ 

\subsection{Weak formulation}

To derive the weak formulation of (\ref{eq:bvp1}), we introduce the
classical Sobolev space 
\begin{equation}
H^{1}(D)=\left\{ u\in L^{2}(D):\nabla u\in\left(L^{2}(D)\right)^{2}\right\} 
\end{equation}
and define the spaces of trial and test functions as,
\begin{eqnarray}
\boldsymbol{V} & = & \left\{ \boldsymbol{v}\in(H^{1}(D))^{2}:\boldsymbol{v}=\boldsymbol{u}_{D}\ \text{on }\Gamma_{D}\right\} \label{eq:wf1}\\
\boldsymbol{V}_{0} & = & \left\{ \boldsymbol{v}\in(H^{1}(D))^{2}:\boldsymbol{v}=\boldsymbol{0}\ \text{on }\Gamma_{D}\right\} \label{eq:wf2}
\end{eqnarray}
The weak formulation of (\ref{eq:bvp1}) is obtained by multiplying
the strong form with an arbitrary test function $\boldsymbol{v}\in \boldsymbol{V}_{0}$
and integrating over $D$, followed by an application of the divergence
theorem and insertion of the Neumann boundary condition. Accordingly,
we arrive at the following weak form.

\subsubsection{Continuous Weak Formulation }

Find $\boldsymbol{u}\in V$ such that 
\begin{eqnarray}
a(\boldsymbol{u};\boldsymbol{v}) & = & l(\boldsymbol{v})\quad\text{for all }\boldsymbol{v}\in \boldsymbol{V}_{0},\label{eq:wf3}
\end{eqnarray}
here the semi-linear form $a(\cdot;\cdot)$ and the linear functional
$l(\cdot)$ are defined by
\begin{eqnarray}
a(\boldsymbol{u};\boldsymbol{v}) & = & \int_{D}\frac{\mathbb{E}[\boldsymbol{\epsilon\left(u\right)}]}{\left\{ 1-\left(\beta\left|\mathbb{E}^{1/2}[\boldsymbol{\boldsymbol{\epsilon\left(u\right)}}]\right|\right)^{\alpha}\right\} ^{1/\alpha}}:\boldsymbol{\varepsilon}(\boldsymbol{v})\,\mathrm{d}\boldsymbol{x},\label{eq:wf4}\\
l(\boldsymbol{v}) & = & \int_{D}\boldsymbol{f}\cdot\boldsymbol{v}\,\mathrm{d}\boldsymbol{x}+\int_{\Gamma_{N}}\boldsymbol{g}\cdot\boldsymbol{v}\,\mathrm{d}s\label{eq:wf5}
\end{eqnarray}

\subsubsection{Well-posedness}

In order for the boundary value problem in (\ref{eq:bvp1}) to be
mathematically well-posed we impose the following set of assumptions,
each of which has a clear physical interpretation. The analytical
framework adopted here is largely inspired by the functional-analytic
setting developed by Beck et al. \cite{beck2017existence} for implicit
constitutive models.
\begin{enumerate}
\item \textbf{Material homogeneity and model parameters.} We assume that
the parameters of the constitutive law are spatially uniform and physically
admissible. In particular, the nonlinear strain-limiting parameters
$\alpha$ and $\beta$ are taken as constants throughout the body
and the elastic moduli $\lambda$ and $\mu$ are also positive constants,
corresponding to a homogeneous material. Moreover, in the limiting
case $\beta\rightarrow0^{+}$ the nonlinear model is required to converge
to the classical linear elasticity model, providing a consistency
check between the strain-limiting and linear theories.
\item \textbf{Equilibrium condition}: When the problem is posed purely in
terms of tractions, i.e., when the Dirichlet boundary is empty, the
externally applied loads must satisfy the global balance condition.
Physically, this means that the net force on the body must vanish
so that the system can be in static equilibrium. Mathematically, this
is expressed by the relation,
\begin{eqnarray}
\int_{D}\boldsymbol{f}d\boldsymbol{\boldsymbol{x}}+\int_{\partial D}\boldsymbol{g}d\boldsymbol{\boldsymbol{s}} & = & \boldsymbol{0.}\label{eq:wf6}
\end{eqnarray}
\item \textbf{Boundary Data Regularity: }The prescribed displacement field
on the Dirichlet boundary, $\boldsymbol{u}_{D}$ is assumed to be
sufficiently regular so that it can be extended into the interior
of the domain. Concretely, $\boldsymbol{u_{D}}\in\left(W^{1,1}\left(D\right)\right)^{2}$
and also \foreignlanguage{american}{$\boldsymbol{\epsilon}\left(\boldsymbol{u_{D}}\left(\boldsymbol{x}\right)\right)$
is contained in a compact set in $\mathbb{R}^{2\times2}$.}
\end{enumerate}
\selectlanguage{american}%
\textit{Theorem : }Let us consider a bounded, connected, Lipschitz
domain\textit{ $D\in\mathbb{R}^{2}$ }which consists of the\textit{
}Dirichlet boundary $\Gamma_{D}$ and the Neumann boundary $\Gamma_{N}$
such that, $\Gamma_{D}\cap\Gamma_{N}=\emptyset$ and $\overline{\Gamma_{D}\cup\Gamma_{N}}=\partial D$.
Given the vector field $\boldsymbol{f}:D\rightarrow\mathbb{R}^{2}$,
a traction vector $\boldsymbol{g}:\Gamma_{N}\rightarrow\mathbb{R}^{2}$
and displacement $\boldsymbol{u_{D}}:\Gamma_{D}\rightarrow\mathbb{R}^{2}$
with the mapping $\mathcal{\bar{F}}:Sym\left(\mathbb{R}^{2\times2}\right)\rightarrow Sym\left(\mathbb{R}^{2\times2}\right)$
for the boundary value problem ((\ref{eq:bvp1})) with the assumptions
stated earlier, then there exists a pair $\left(\mathbf{u,\boldsymbol{\sigma}}\right)\in\left(W^{1,1}\left(D\right)\right)^{2}\times Sym\left(L^{1}\left(D\right)^{2\times2}\right)$
satisfying the weak form,

\begin{eqnarray*}
\int_{D}\boldsymbol{\sigma}(\boldsymbol{\epsilon}(\boldsymbol{u})):\boldsymbol{\epsilon}(\boldsymbol{w})\,d\boldsymbol{x} & = & \int_{D}\boldsymbol{f}\cdot\boldsymbol{w}\,d\boldsymbol{x}+\int_{\Gamma_{N}}\boldsymbol{g}\cdot\boldsymbol{w}\,d\boldsymbol{s},\quad\forall\boldsymbol{w}\in \boldsymbol{V}_{0}.
\end{eqnarray*}
\bigskip{}
\foreignlanguage{english}{In order to compute approximate solutions,
we now introduce a conforming finite element discretization based
on continuous, piecewise polynomial vector fields.}

\subsubsection{Finite element discretization}

\selectlanguage{english}%
Let $\{\mathcal{T}_{h}\}_{h>0}$ be a family of shape-regular partitions
of the domain $D$, where each mesh $\mathcal{T}_{h}$ consists of
quadrilaterals $K$ such that,
\begin{eqnarray*}
\overline{D}=\bigcup_{K\in\mathcal{T}_{h}}\overline{K}
\end{eqnarray*}
and $h:=\max_{K\in\mathcal{T}_{h}}h_{K}$ denotes the global mesh
size with $h_{K}$ the diameter of the element $K$. The crack $\Gamma_{c}$
is treated as an internal traction-free boundary. On each element
$K\in\mathcal{T}_{h}$ we work with vector-valued polynomial shape
functions in $\left(\mathbb{Q}_{1}(K)\right)^{2}$. The discrete trial
and test spaces are then defined as,
\begin{eqnarray*}
\boldsymbol{V}_{h} & = & \bigl\{\boldsymbol{v}_{h}\in C^{0}(\overline{D})^{2}:\boldsymbol{v}_{h}|_{K}\in\left(\mathbb{Q}_{1}(K)\right)^{2}\ \forall K\in\mathcal{T}_{h},\ \boldsymbol{v}_{h}=\boldsymbol{u}_{D}\text{ on }\Gamma_{D}\bigr\},\\
\boldsymbol{V}_{h,0} & = & \bigl\{\boldsymbol{v}_{h}\in C^{0}(\overline{D})^{2}:\boldsymbol{v}_{h}|_{K}\in\left(\mathbb{Q}_{1}(K)\right)^{2}\ \forall K\in\mathcal{T}_{h},\ \boldsymbol{v}_{h}=\boldsymbol{0}\text{ on }\Gamma_{D}\bigr\},
\end{eqnarray*}

\paragraph{Discrete weak form.}

Let $a(\cdot;\cdot)$ and $l(\cdot)$ denote the semi-linear form
and the linear functional introduced in(\ref{eq:wf4}) and (\ref{eq:wf5}).
The finite element approximation of the boundary value problem consists
in finding a discrete displacement field $\boldsymbol{u}_{h}\in \boldsymbol{V}_{h}$
such that,
\begin{eqnarray}
a(\boldsymbol{u}_{h};\boldsymbol{v}_{h})=l(\boldsymbol{v}_{h}) & \qquad\text{for all } & \boldsymbol{v}_{h}\in  \boldsymbol{V}_{h,0}.\label{eq:dw1}
\end{eqnarray}
In expanded form, this relation can be written as,
\begin{eqnarray}
\int_{D}\frac{\mathbb{E}[\boldsymbol{\epsilon}\left(\boldsymbol{u}_{h}\right)]}{\left\{ 1-\left(\beta\left|\mathbb{E}^{1/2}[\boldsymbol{\epsilon}\left(\boldsymbol{u}_{h}\right)]\right|\right)^{\alpha}\right\} ^{1/\alpha}}:\boldsymbol{\varepsilon}(\boldsymbol{v}_{h})\,\mathrm{d}\boldsymbol{x} & = & \int_{D}\boldsymbol{f}\cdot\boldsymbol{v}_{h}\,\mathrm{d}x+\int_{\Gamma_{N}}\boldsymbol{g}\cdot\boldsymbol{v}_{h}\,\mathrm{d}s,\nonumber \\
 &  & \forall\boldsymbol{v}_{h}\in V_{h,0}.\label{eq:dw2}
\end{eqnarray}

\paragraph{Picard linearization.} The formulation given in (\ref{eq:dw2})
It is non-linear due to the term on the left-hand side. So we implement Picard's
iterative method to deal with this non-linearity. Starting from an
initial guess $\boldsymbol{u}_{h}^{(0)}$ (solution of the linear
elastic model corresponding to $\beta=0$), we generate a sequence
of approximations $\{\boldsymbol{u}_{h}^{(k)}\}_{k\ge0}\subset V_{h}$
by solving a sequence of linearized variational problems. Given an
iterate $\boldsymbol{u}_{h}^{(k)}\in V_{h}$, we define the scalar
coefficient,
\begin{eqnarray}
\xi^{(k)}(x) & = & \frac{1}{\bigl(1-(\beta\|\mathbb{E}^{1/2}[\boldsymbol{\boldsymbol{\epsilon}}(\boldsymbol{u}_{h}^{(k)})]\|)^{\alpha}\bigr)^{1/\alpha}},\qquad x\in D,\label{eq:pl1}
\end{eqnarray}
The next iterate $\boldsymbol{u}_{h}^{(k+1)}\in V_{h}$ is then defined
as the solution of the linear problem,

\begin{eqnarray}
\int_{D}\xi^{(k)}(x)\,\mathbb{E}\big[\boldsymbol{\boldsymbol{\epsilon}}(\boldsymbol{u}_{h}^{(k+1)})\big]:\boldsymbol{\boldsymbol{\epsilon}}(\boldsymbol{v}_{h})\,\mathrm{d}\boldsymbol{x} & = & \int_{D}\boldsymbol{f}\cdot\boldsymbol{v}_{h}\,\mathrm{d}\boldsymbol{x}+\int_{\Gamma_{N}}\boldsymbol{g}\cdot\boldsymbol{v}_{h}\,\mathrm{d}\boldsymbol{s},\label{eq:pl2}\\
 &  & \quad\forall\boldsymbol{v}_{h}\in V_{h,0}.\nonumber 
\end{eqnarray}
The iteration is terminated once a chosen convergence criterion is
met, for example when the residual falls below a prescribed tolerance.
\selectlanguage{american}%

\subsection{Geometry and boundary conditions}

\selectlanguage{english}%
The configuration under consideration is a rectangular plate of unit
dimensions in the $(x,y)$-plane. Without loss of generality, the
reference domain is taken to be 
\[
D=[0,1]\times[0,1]\subset\mathbb{R}^{2},
\]
and the plate contains a pre-existing edge crack. The crack is modeled
as a line segment $\Gamma_{c}$ extending inward from the right boundary
along the mid-line, 
\[
\Gamma_{c}=\{(x,y)\in\overline{D}:0.5\leq x\leq1,\;y=0.5\},
\]
and is treated as an internal traction-free boundary. The outer boundary
$\partial D$ is decomposed into four segments, $\partial D=\Gamma_{1}\cup\Gamma_{2}\cup\Gamma_{3}\cup\Gamma_{4}\cup\Gamma_{c}$
The top and bottom edges $\Gamma_{D}:=\Gamma_{1}\cup\Gamma_{3},$
form the Dirichlet boundary on which vertical displacements are prescribed,
while the remaining edges together with the crack boundary,
\begin{eqnarray*}
\Gamma_{N} & = & \Gamma_{2}\cup\Gamma_{4}\cup\Gamma_{c},
\end{eqnarray*}
are treated as traction-free Neumann boundaries.
\begin{figure}[H]
\centering
\includegraphics[scale=0.5]{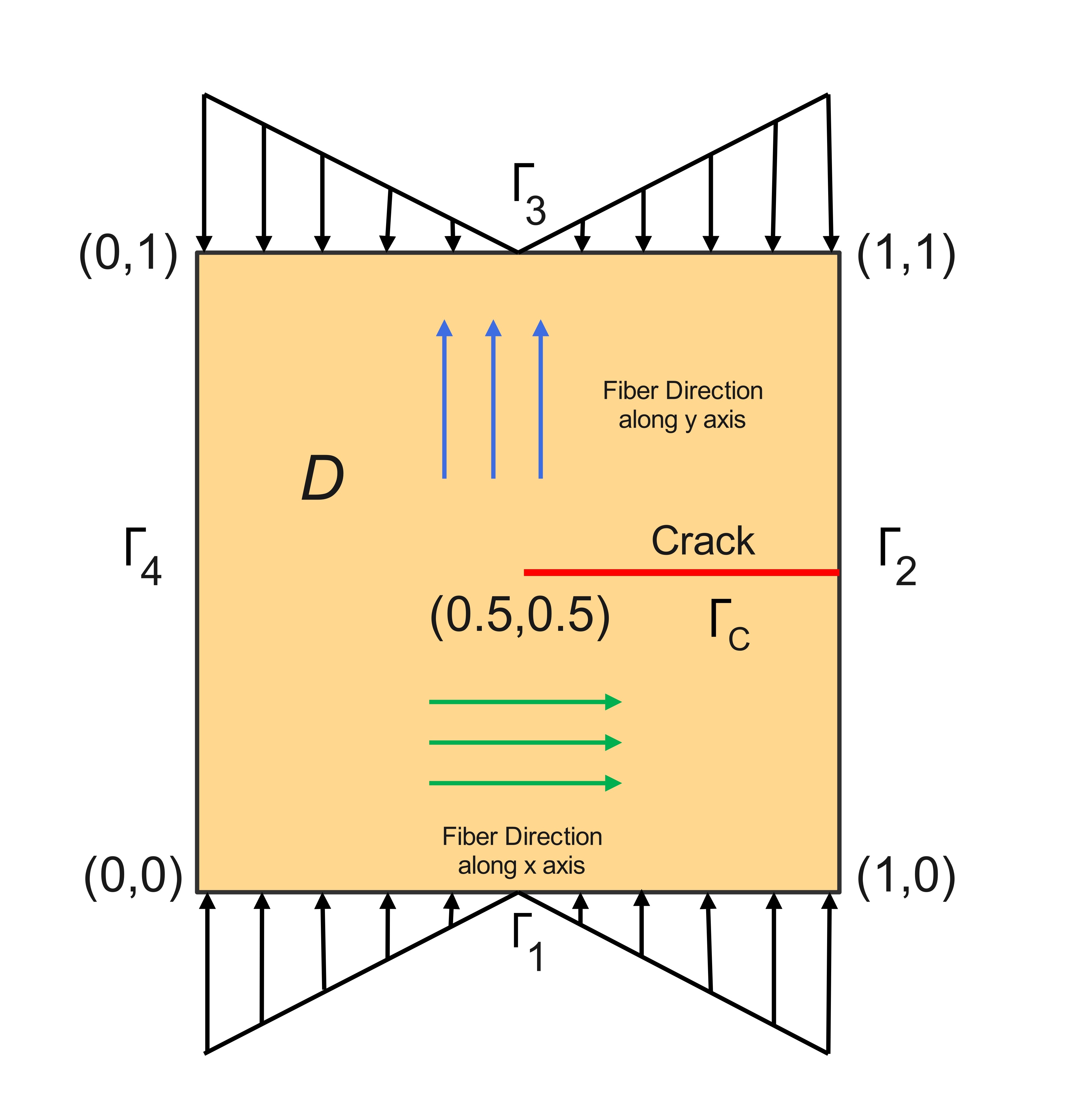}
\caption{Schematic representation of the problem geometry and boundary conditions. A rectangular domain with an edge crack is subjected to piecewise linear vertical displacements on the top and bottom boundaries, with traction-free conditions on the remaining edges.}
\end{figure}

The key distinguishing feature of the present study is the piecewise
slope loading prescribed on the top and bottom edges of the cracked
plate. We prescribe purely vertical displacements $u_{2}$ on these
edges and the displacement component is chosen to vary piecewise linearly
along $\Gamma_{1}$ and $\Gamma_{3}$ with a change of slope at the
mid point, $x=0.5$ while the horizontal displacement component $u_{1}$
is left unconstrained. Introducing a slope parameter $\ell>0$, the
boundary conditions for the top and bottom boundary are as follows.
\begin{eqnarray*}
u_{2}\left(x,1\right) & = & \begin{cases}
\ell\,(x-0.5), & 0\le x\le0.5,\\[2pt]
\ell\,(0.5-x), & 0.5\le x\le1,
\end{cases}\text{on }\Gamma_{3},\\
u_{2}\left(x,0\right) & = & \begin{cases}
\ell\,(0.5-x), & 0\le x\le0.5,\\[2pt]
\ell\,(x-0.5), & 0.5\le x\le1,
\end{cases}\text{on }\Gamma_{1},
\end{eqnarray*}
The left $\left(\Gamma_{4}\right)$ and right $\left(\Gamma_{2}\right)$
boundaries together with the crack boundary $\Gamma_{C}$ are treated
as the part of the Neumann boundary $\Gamma_{N}$ and considered as
traction-free,
\begin{eqnarray*}
\boldsymbol{\sigma n} & = & \boldsymbol{0}\quad\text{on }\Gamma_{2}\cup\Gamma_{4}\cup\Gamma_{c},
\end{eqnarray*}
The only externally imposed deformation is therefore generated by
the boundary conditions on $\Gamma_{1}$ and $\Gamma_{3}$ by the
vertical displacement component.

\section{Numerical Results and Discussion}

In this section, we present numerical experiments designed to assess
the effect of the piecewise slope boundary loading on the crack-tip
fields in a transversely isotropic strain-limiting body. The computational
experiments in this study are carried out using the open-source finite
element library \foreignlanguage{american}{deal.II \cite{2023dealii,dealii2019design}.
The domain is discretized with a bilinear finite element space. Convergence
is monitored using the residual associated with the weak form of the
governing equations.} It is computed using the following formula,

\begin{eqnarray}
\mathtt{R}\left(\boldsymbol{u}_{h}^{(k)}\right) & = & \int_{D}\frac{\mathbb{E}[\boldsymbol{\epsilon}\left(\boldsymbol{u}_{h}^{(k)}\right)]}{\left\{ 1-\left(\beta\left|\mathbb{E}^{1/2}[\boldsymbol{\epsilon}\left(\boldsymbol{u}_{h}^{(k)}\right)]\right|\right)^{\alpha}\right\} ^{1/\alpha}}:\boldsymbol{\varepsilon}(\boldsymbol{v}_{h})\,\mathrm{d}\boldsymbol{x}\label{eq:res1}
\end{eqnarray}
We compute stress, strain, and the strain-energy density near the crack
tip. In what follows, the figures correspond to the numerical plots
obtained from the simulations.
Algorithm 1 summarizes the numerical procedure employed in this work.
The solution is obtained through a Picard iteration on each mesh,
while adaptive refinement driven by a Kelly error estimator progressively
improves resolution near the crack tip and regions of high deformation
induced by the piecewise-slope boundary loading.\medskip{}

\begin{algorithm}[H]
    \caption*{\textbf{Algorithm:} Adaptive FEM for Strain-Limiting Crack Model under Piecewise-Slope Loading}
    
    \begin{algorithmic}[1]
        % Using tabular to align inputs cleanly under a single "Input:" label
        \Require 
        \begin{tabular}[t]{@{}l@{}}
            Material parameters $(\lambda, \mu, \bar{\gamma}, \alpha, \beta)$ and slope parameter $\ell > 0$ \\
            Tolerance $\varepsilon_{\mathrm{tol}}$, Maximum Picard iterations $k_{\max}$ \\
            Number of adaptive refinement cycles $N$
        \end{tabular}
        
        \Statex % Vertical spacing
        \State \textbf{Step 1: Initialization}
        \State Generate initial mesh $\mathcal{T}_{0}$ and tag boundary $\Gamma_{c}$
        \State Construct finite element space $V_{h}$ on $\mathcal{T}_{0}$
        \State Compute initial guess $u_{h}^{(0)}$ using linear elasticity (set $\beta=0$)
        
        \Statex
        \For{$m = 0, \dots, N$}
            \If{$m > 0$}
                \State \textbf{Step 2: Solution Transfer}
                \State Interpolate solution $u_{h}^{(m-1)}$ onto the current mesh $\mathcal{T}_{m}$
            \EndIf
            
            \Statex
            \State \textbf{Step 3: Picard Iteration}
            \For{$k = 0, \dots, k_{\max}$}
                \State Assemble system matrix $A(\xi^{(k)})$ and RHS $f$
                \State Solve linear system: $A(\xi^{(k)}) \, u_{h}^{(k+1)} = f$
                \State Compute residual norm $R = \|\mathtt{R}(u_{h}^{(k+1)})\|$
                \If{$R < \varepsilon_{\mathrm{tol}}$}
                    \State \textbf{break} \Comment{Convergence achieved}
                \EndIf
            \EndFor
            
            \Statex
            \State \textbf{Step 4: Adaptation}
            \If{$m < N$}
                \State Compute error indicators $\eta_{K}$ for each element $K$ (Kelly estimator)
                \State Mark elements for refinement (e.g., top $30\%$ strategy)
                \State Refine mesh to generate $\mathcal{T}_{m+1}$
            \EndIf
        \EndFor
        
        \Statex
        \State \textbf{Step 5: Post-processing}
        \State Compute stress, strain, and energy density fields on the final mesh
    \end{algorithmic}
\end{algorithm}

\subsection{Mesh Convergence Analysis}

A reliable finite element approximation of fracture fields requires
not only accurate resolution of localized deformation near the crack
tip but also demonstrable convergence of the computed quantities of
interest as the mesh is refined. The adaptive refinement strategy
follows the standard strategy: 
\[
SOLVE \rightarrow  ESTIMATE \rightarrow
MARK \rightarrow REFINE.
\]
\begin{enumerate}
\item \textbf{Solve:} For a given mesh $\mathcal{T}_{h}$, compute the discrete
solution $u_{h}$.
\item \textbf{Estimate:} Evaluate an a posteriori error indicator $\eta_{K}$
on each element $K\in\mathcal{T}_{h}$. In the present work we use
a Kelly error estimator (jump-based) applied to the displacement solution;
the indicator is concentrated near the crack tip and along the crack
faces.
\item \textbf{Mark:} Select a subset of elements for refinement using a
fixed-fraction strategy: elements with the largest indicators are
marked so that approximately a prescribed fraction of the total error
is accounted for.
\item \textbf{Refine :} Refine the marked elements (those with the largest
indicators) in regions away from the crack tip. The discrete solution
is transferred to the new mesh by interpolation.
\end{enumerate}
For a given mesh $\mathcal{T}_{h}$, the boundary value problem is
first solved using the Picard linearization described in Section~3.2.
Element-wise a posteriori error indicators $\eta_{K}$ are then computed
using the residual-based Kelly estimator implemented in \texttt{deal.II}.
These indicators are typically largest near the crack tip and along
the crack faces, but under the present piecewise-slope boundary conditions
they also detect regions of high curvature near the slope-reversal
points on the top and bottom boundaries. Elements contributing approximately
the top 30\% of the total estimated error are marked for refinement,
and the mesh is locally refined without coarsening. The discrete solution
is subsequently transferred to the new mesh by interpolation, and
the procedure is repeated. As a representative test case, we consider material fibers aligned
with the $x$ - axis and fix the model parameters to $\alpha=0.5$
and $\beta=0.05$. The computation is initiated on a relatively coarse
mesh with $256$ active cells and $594$ degrees of freedom. A sequence
of $12$ adaptive refinement cycles is then performed (cycles $k=0,\dots,11$).
After the final adaptive cycle the mesh contains $25147$ active cells
and $52346$ degrees of freedom, representing over an order of magnitude
increase in resolution compared with the initial grid. The below table
summarizes the evolution of the mesh complexity with refinement.
\medskip{}

\begin{table}[H]
\centering
\centering \caption{Adaptive refinement for $\alpha=0.5$ and $\beta=0.05$, including
the monitored displacement value $u_{2}(0.75,0.5)$ at each cycle.}
\label{tab:mesh_info} %
\begin{tabular}{cccc}
\hline 
Cycle $i$ & Active cells & DoFs & $u_{2}(0.75,0.5)$\tabularnewline
\hline 
0 & 256 & 594 & $-0.11417338944$\tabularnewline
1 & 343 & 818 & $-0.11445971655$\tabularnewline
2 & 508 & 1212 & $-0.11486387531$\tabularnewline
3 & 778 & 1826 & $-0.11507118568$\tabularnewline
4 & 1210 & 2794 & $-0.11524236650$\tabularnewline
5 & 1879 & 4230 & $-0.11524689579$\tabularnewline
6 & 2887 & 6394 & $-0.11530167631$\tabularnewline
7 & 4456 & 9712 & $-0.11532721825$\tabularnewline
8 & 6949 & 14920 & $-0.11533615721$\tabularnewline
9 & 10612 & 22554 & $-0.11535055907$\tabularnewline
10 & 16204 & 34056 & $-0.11535906484$\tabularnewline
11 & 25147 & 52346 & $-0.11536127909$\tabularnewline
\hline 
\end{tabular}
\end{table}

The monotonic stabilization of $u_{2}(0.75,0.5)$ across refinement
cycles provides a simple yet effective scalar indicator of convergence.
To quantify this more systematically, the displacement obtained on
the finest mesh (cycle 11) is adopted as a reference value. The relative
error at cycle $i$ is defined as 
\[
\boldsymbol{\epsilon}_{i}=\frac{\left|u_{2}^{(i)}(0.75,0.5)-u_{2}^{\mathrm{ref}}(0.75,0.5)\right|}{\left|u_{2}^{\mathrm{ref}}(0.75,0.5)\right|}.
\]
\begin{figure}[H]
\centering
\includegraphics[scale=0.2]{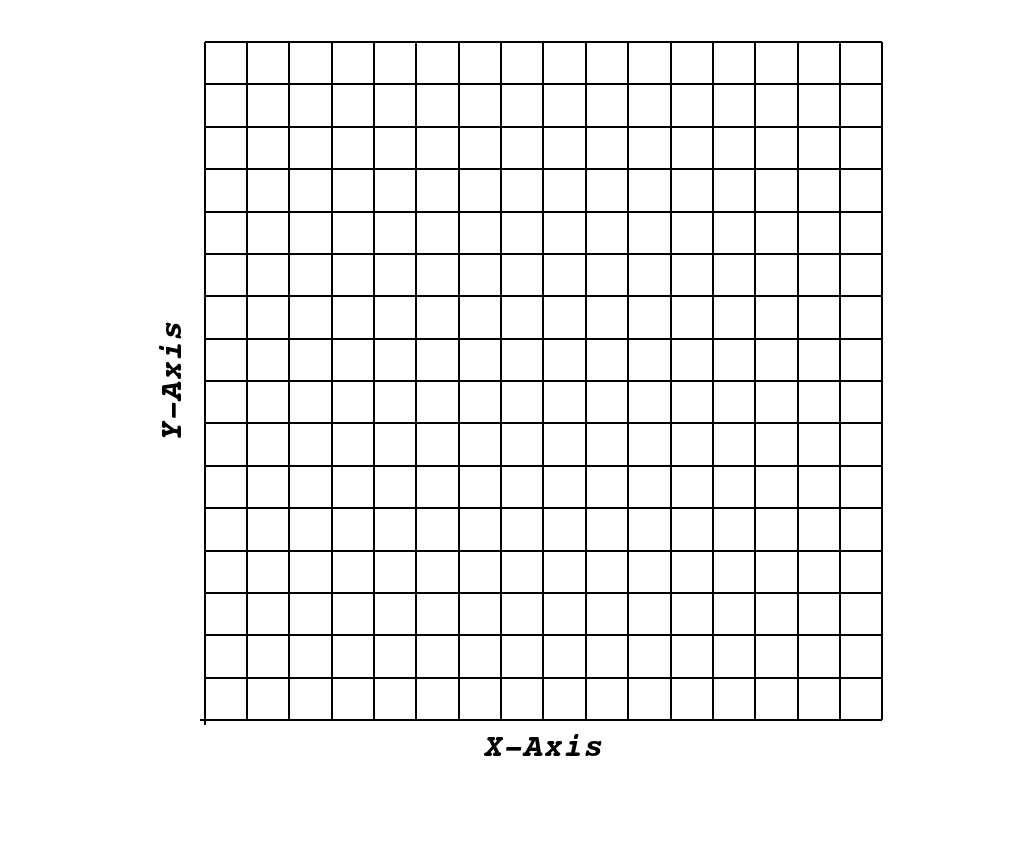}\includegraphics[scale=0.2]{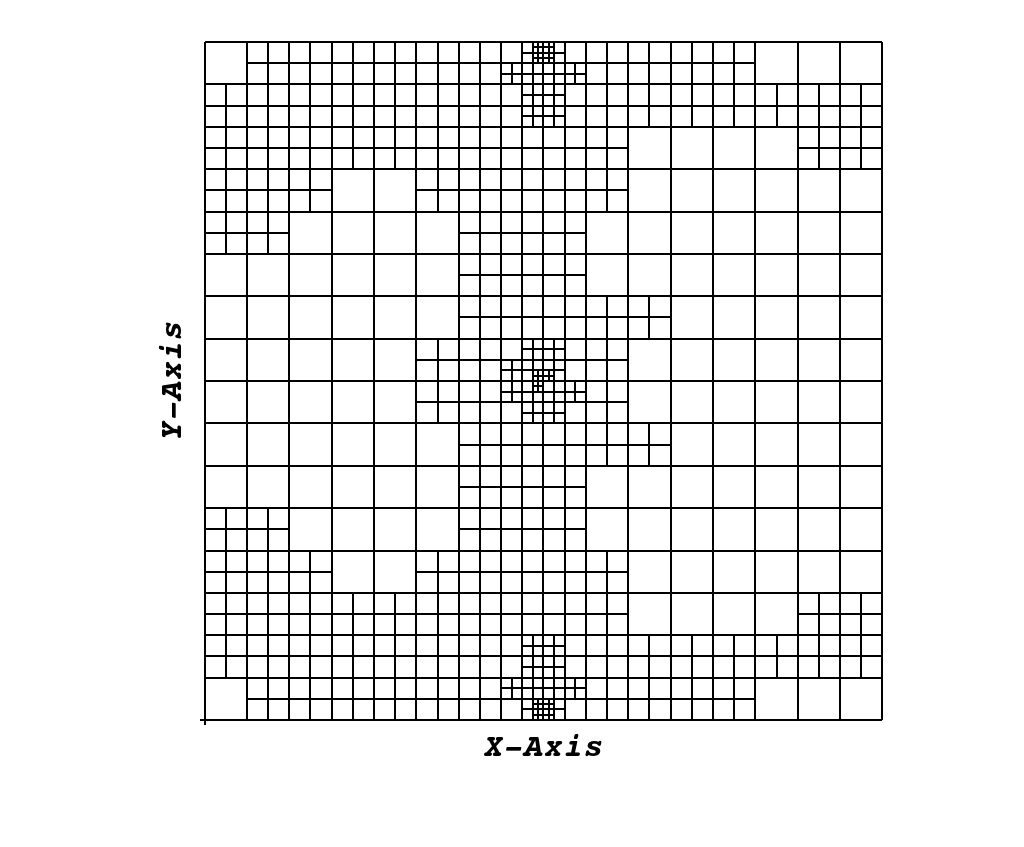}
\caption{Evolution of the adaptive mesh refinement. Left: The initial coarse mesh (Cycle 0). Right: The refined mesh at Cycle 3, showing initial localization of elements around the crack tip and regions of boundary slope discontinuities.}
\end{figure}
\bigskip{}
\begin{figure}[H]
\centering
\includegraphics[scale=0.2]{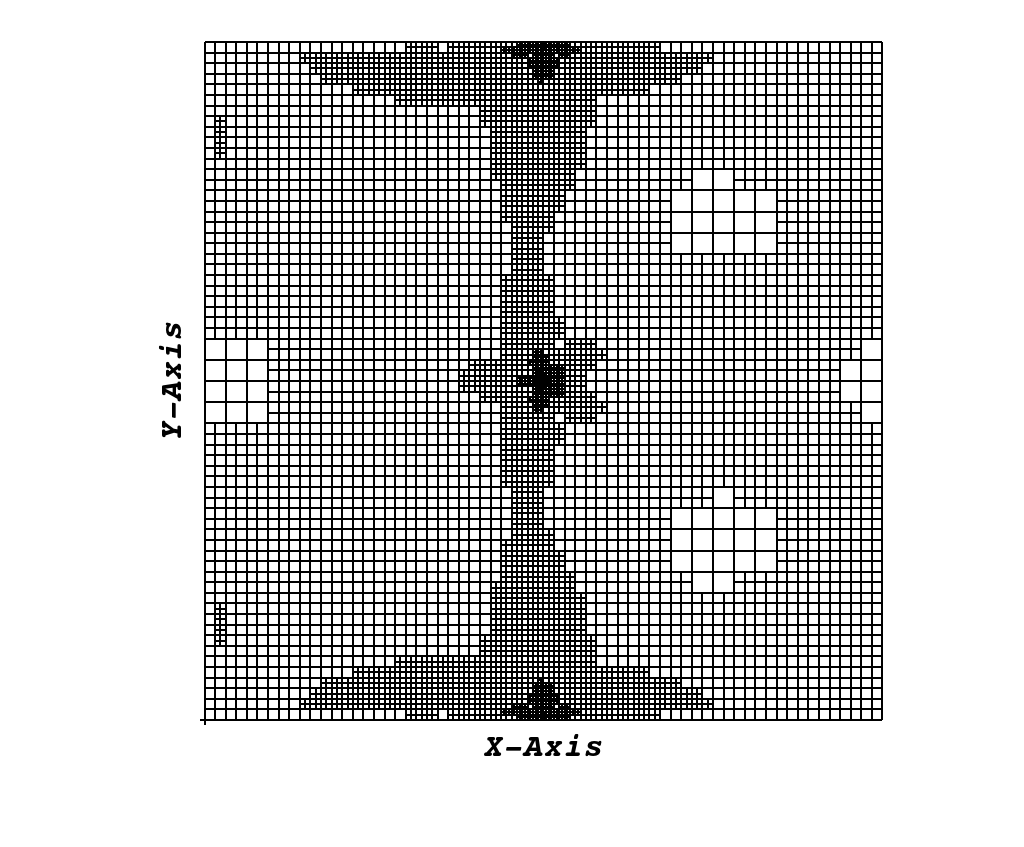}\includegraphics[scale=0.2]{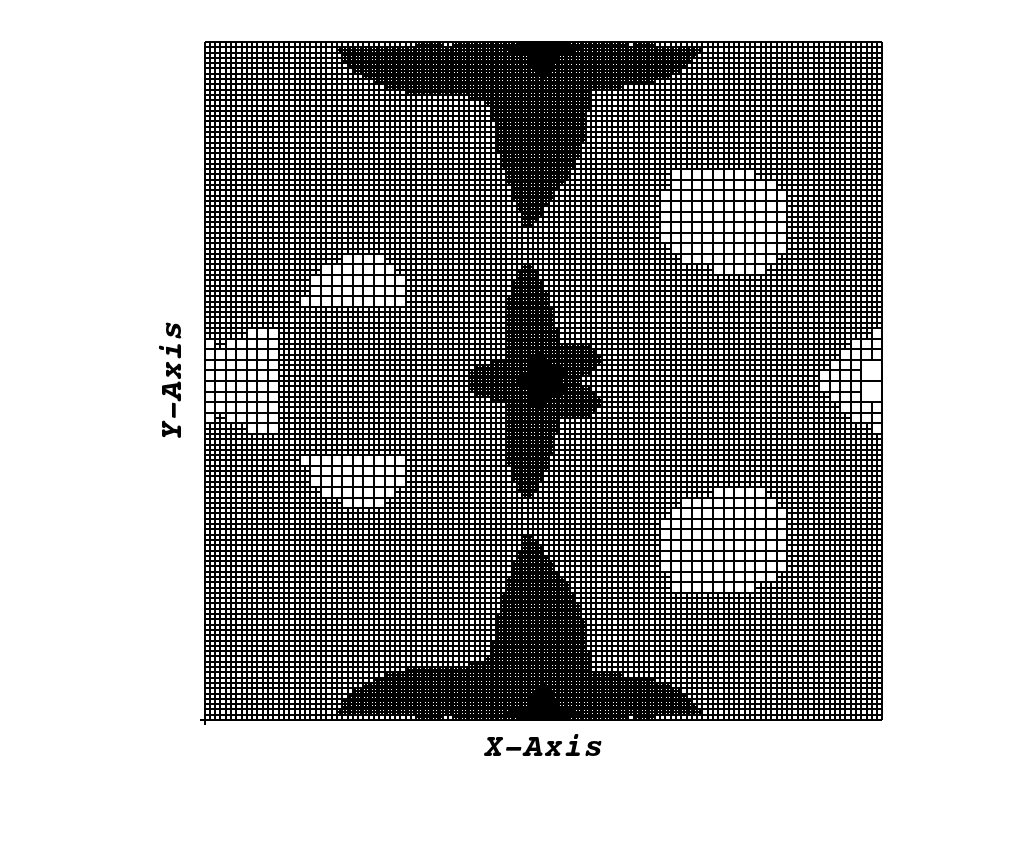}
\caption{Advanced stages of adaptive mesh refinement. Left: The mesh at Cycle 8. Right: The finest mesh at Cycle 11. Note the highly focused refinement along the crack line and the emergence of specific refinement patterns extending from the top and bottom boundaries, driven by the piecewise slope loading.}
\end{figure}

A distinctive feature of the present problem, compared with more conventional
uniform or smoothly varying loads, is the presence of slope discontinuities
in the prescribed boundary displacement on $\Gamma_{1}$ and $\Gamma_{3}$.
The change of sign in the slope at $x=0.5$ gives rise to pronounced
variations in the displacement gradient near the midpoints of the
top and bottom boundaries. The adaptive refinement indicators are
therefore influenced not only by the crack geometry but also by the
non-smooth character of the displacement loading. Throughout the refinement
process we consistently observe that cells are marked in three principal
regions: (i) a narrow band surrounding the crack tip, (ii) localized
refinement layers originating from the slope-reversal points on the
top and bottom boundaries, and (iii) the transition region in which
these boundary layers interact with the crack-tip field. This behavior
is in sharp contrast to the purely uniform-loading case, where refinement
is largely confined to a neighborhood of the crack and decays rapidly
away from the tip. Under the piecewise-slope boundary conditions,
the adaptive algorithm automatically detects and resolves the additional
curvature in the displacement field induced by the change in slope,
leading to an anisotropic refinement pattern that follows both the
crack line and the regions of slope change. \begin{figure}[H]
\centering
\includegraphics[scale=0.4]{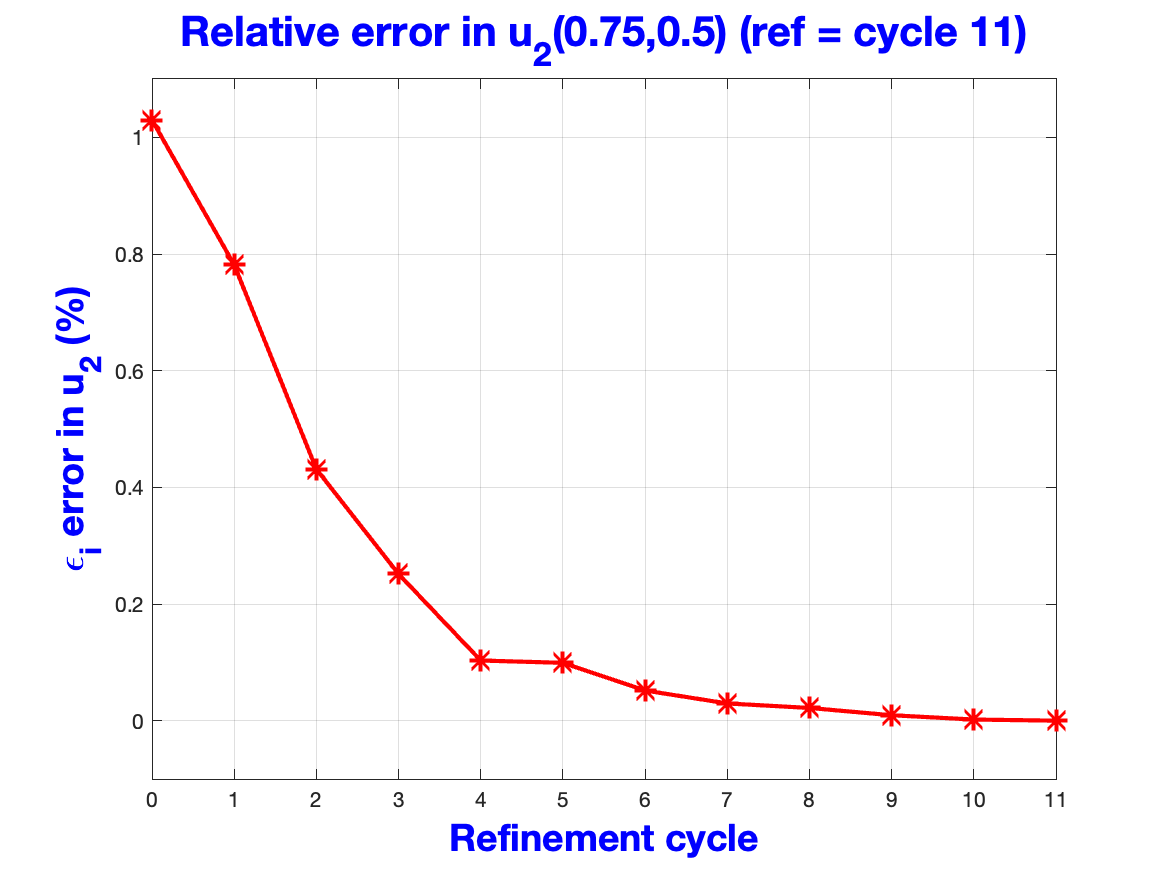}
\caption{Convergence of the vertical displacement component $u_{2}$ at the crack center $(0.75, 0.5)$ as a function of the adaptive refinement cycle. The relative error rapidly decreases and stabilizes below $1\%$ after approximately 6 cycles, indicating mesh independence of the computed solution.\label{fig:error}}
\end{figure}
 Figure~\ref{fig:error} shows the evolution of the relative error
in the vertical displacement component evaluated at the midpoint
of the crack $\left(0.75,0.5\right)$ with respect to the reference
solution obtained on the finest adaptive mesh. A rapid reduction in
the relative error is observed during the early refinement cycles,
followed by a gradual decay as the mesh is further refined. Beyond
approximately cycle 6, the error decreases steadily to below one percent,
indicating that the displacement field has effectively stabilized
at this location. This behavior confirms that the adaptive refinement
strategy efficiently improves solution accuracy not only near singular
features but also in the interior of the domain, leading to mesh-independent
resolution of the displacement field. \begin{center}
\begin{figure}[H]
\centering
\includegraphics[scale=0.25]{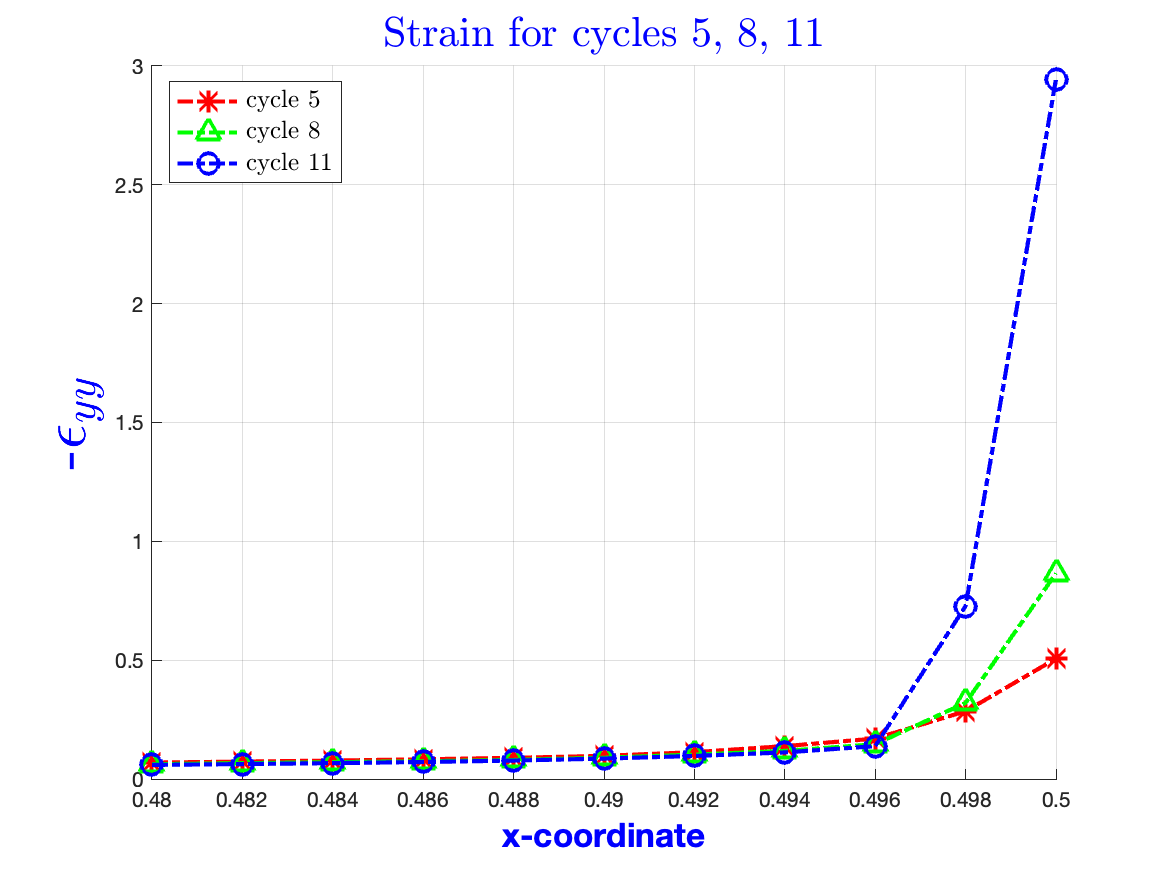}\includegraphics[scale=0.25]{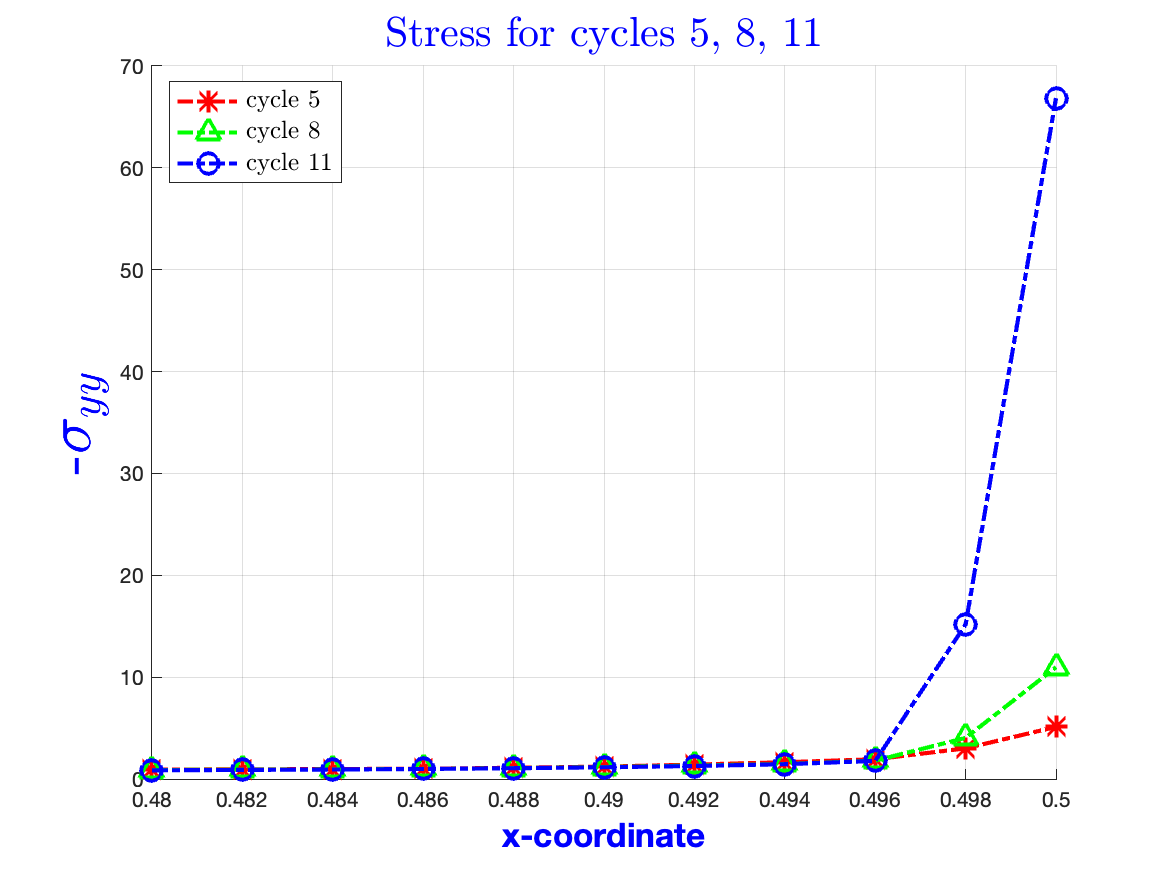}

\caption{Evolution of crack-tip field profiles along the crack line ($x \in [0.48, 0.5]$) for adaptive cycles 5, 8, and 11, with $\beta=0.05$ and $\alpha=0.5$. Left: Vertical strain $\boldsymbol{\epsilon}_{yy}$. Right: Vertical stress $\boldsymbol{\sigma}_{yy}$. The plots illustrate the sharp intensification of stress gradients near the tip as the mesh resolution improves.}

\end{figure}
\bigskip{}
\par\end{center}

\begin{table}[H]
\centering \caption{Crack-tip fields values at $x=0.5$ for $\alpha=0.5$ and $\beta=0.05$
for selected adaptive cycles.}
\label{tab:tip2} 
\centering{}%
\begin{tabular}{cccc}
\hline 
\textbf{Crack tip fields} & \textbf{Cycle 5} & \textbf{Cycle 8} & \textbf{Cycle 11}\tabularnewline
\hline 
$\boldsymbol{\epsilon}_{yy}$ & -0.506503 & -0.865436 & -2.94502\tabularnewline
$\boldsymbol{\sigma}_{yy}$ & -5.11227 & -10.971 & -66.8477\tabularnewline
\hline 
\end{tabular}
\end{table}

From the strain and stress profiles plotted near crack tip $\left(x\in[0.48,\,0.5]\right)$,
both $\boldsymbol{\epsilon}_{yy}$ and $\boldsymbol{\sigma}_{yy}$
remain relatively mild away from the crack tip but rise sharply as
$x\to0.5$. Moreover, this near-tip amplification becomes dramatically
steeper with successive adaptive refinements (cycles $5\rightarrow8\rightarrow11$).
This behavior is precisely what adaptive refinement is designed to
reveal: as the mesh increasingly resolves the crack-tip neighborhood,
the computed fields are able to capture strong local gradients that
are otherwise smeared out on coarser meshes. The growth is particularly
pronounced for the stress field.

\begin{figure}[H]
\centering
\includegraphics[scale=0.4]{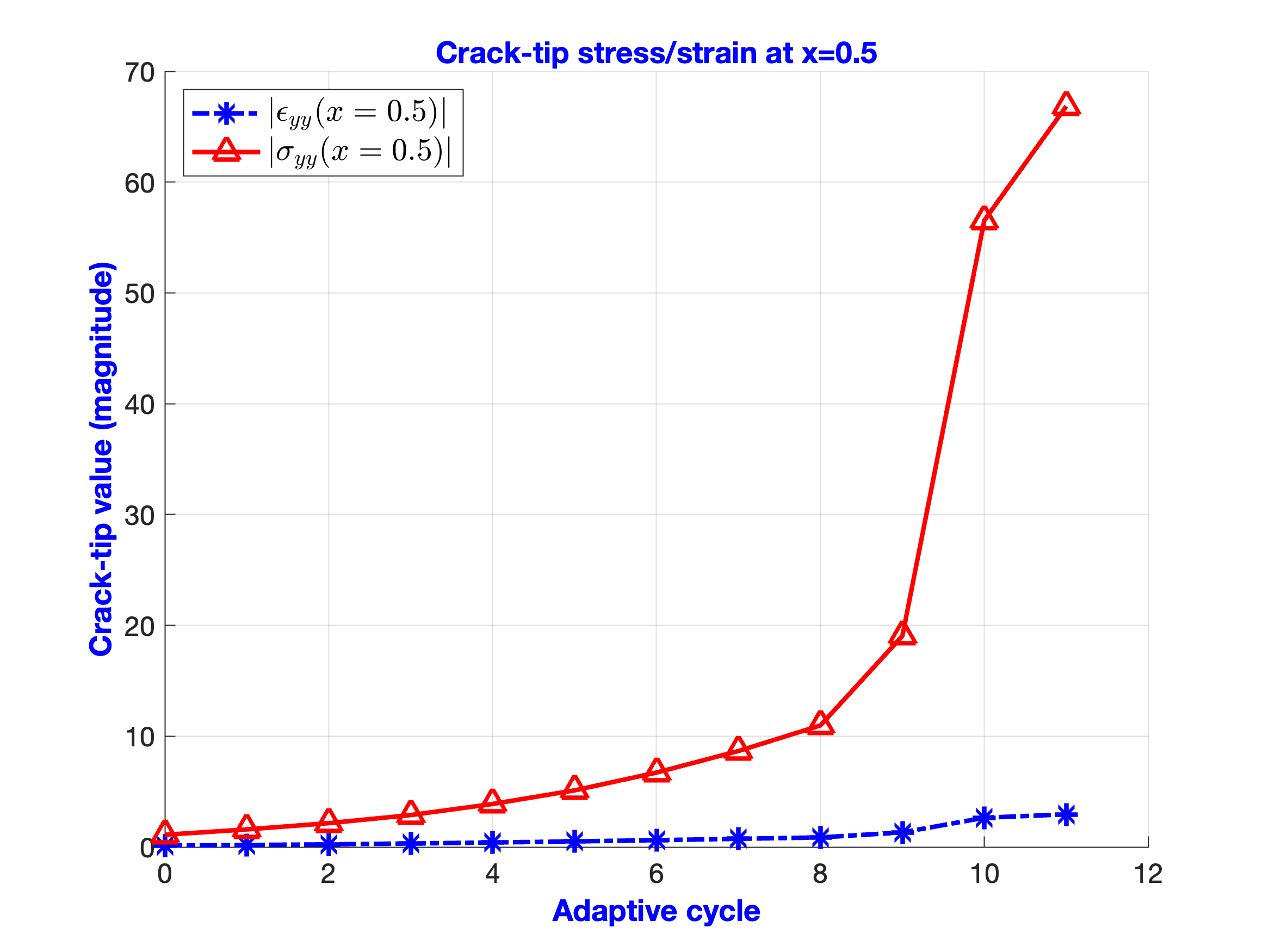}

\caption{Comparison of the growth rates for crack-tip stress and strain across adaptive refinement cycles. The substantial escalation of stress compared to the moderated growth of strain highlights the strain-limiting character of the constitutive model.}
\end{figure}
\bigskip{}
While the crack-tip strain increases from approximately $0.5$ in
cycle~5 to about $2.95$ in cycle~11, the corresponding crack-tip
stress increases from about $5.11$ to nearly $66.85$. The substantially
faster escalation of stress compared to strain is consistent with
the intent of the strain-limiting constitutive model: strains remain
comparatively moderated, whereas stresses can still localize strongly
near the crack tip as the resolution improves. Having established the stability of the numerical scheme and the effectiveness
of the adaptive mesh refinement strategy under piecewise slope boundary
loading, we now investigate the influence of material anisotropy on
the computed crack-tip fields. In particular, we examine how the orientation
of the preferred fiber direction affects the stress, strain, and displacement
distributions in a transversely isotropic strain-limiting elastic
body. We consider two principal fiber orientations:
\begin{itemize}
\item Case $I$: fibers aligned with the crack plane (along the $x$-axis).
\item Case $II$: fibers aligned orthogonal to the crack plane (along the
$y$-axis).
\end{itemize}

\subsection{Case I: Fiber Directions along X axis}

We first consider the case in which the preferred fiber direction
is aligned with the $x$-axis, that is, $\boldsymbol{e_{1}}=(1,0)$
so that the structural tensor is $M=\boldsymbol{e_{1}}\otimes\boldsymbol{e_{1}}$.
This choice corresponds to a material that is stiffer along the horizontal
direction, parallel to the crack plane. \selectlanguage{american}%
\begin{center}

\selectlanguage{english}%
\begin{figure}[H]
\centering
\includegraphics[scale=0.25]{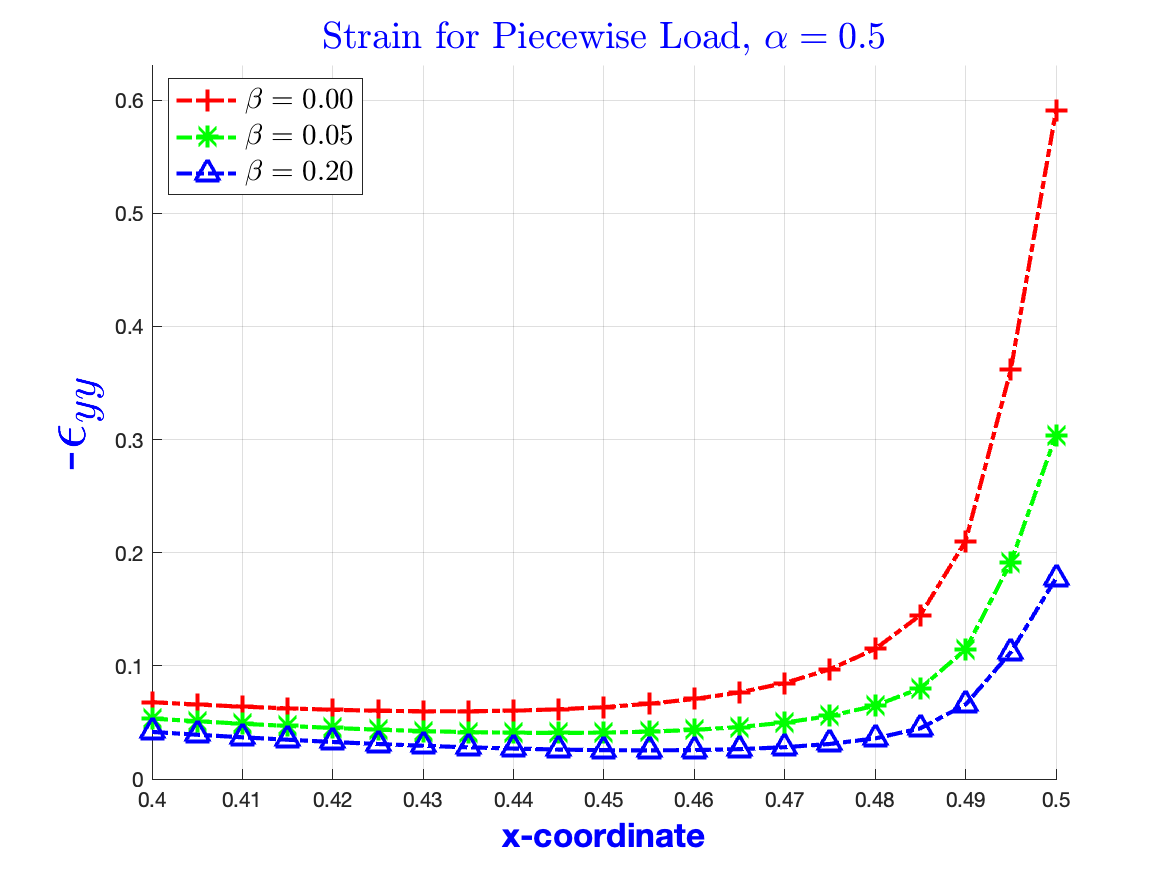}\includegraphics[scale=0.25]{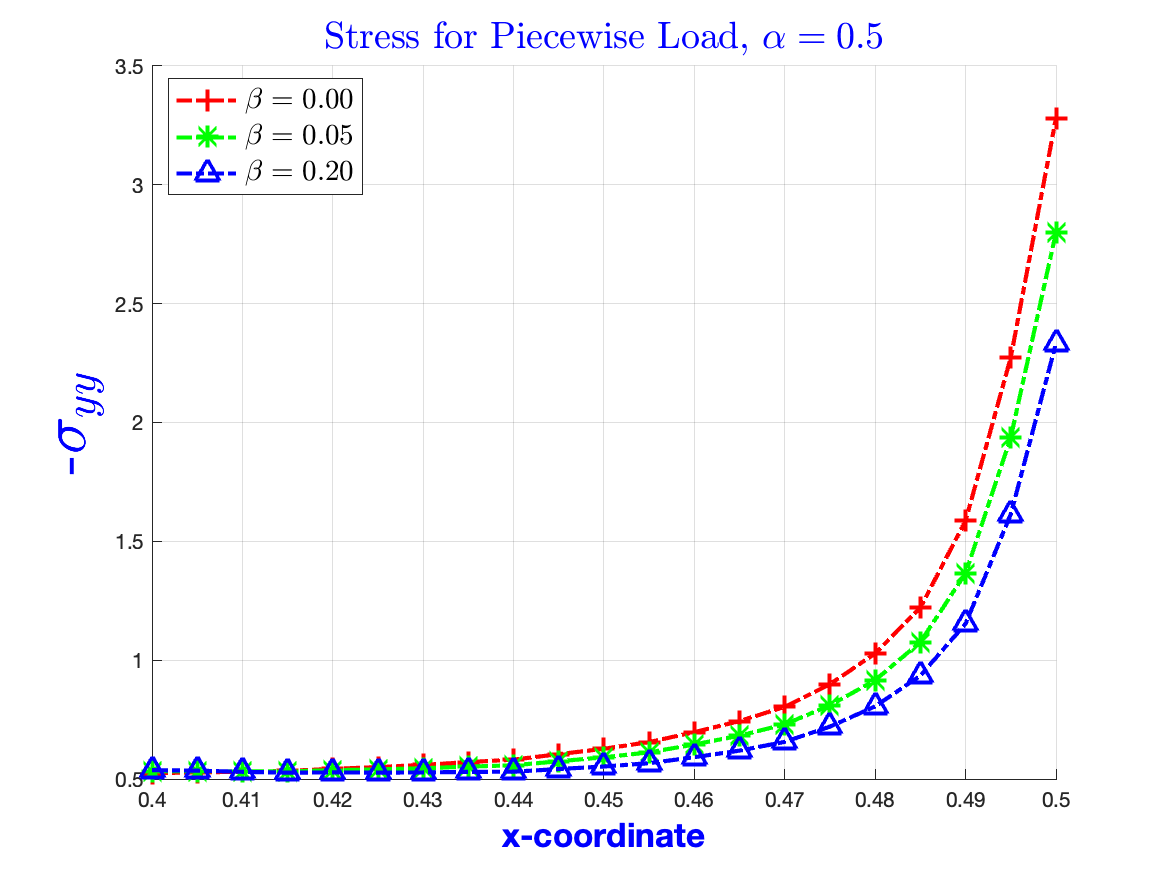}\includegraphics[scale=0.25]{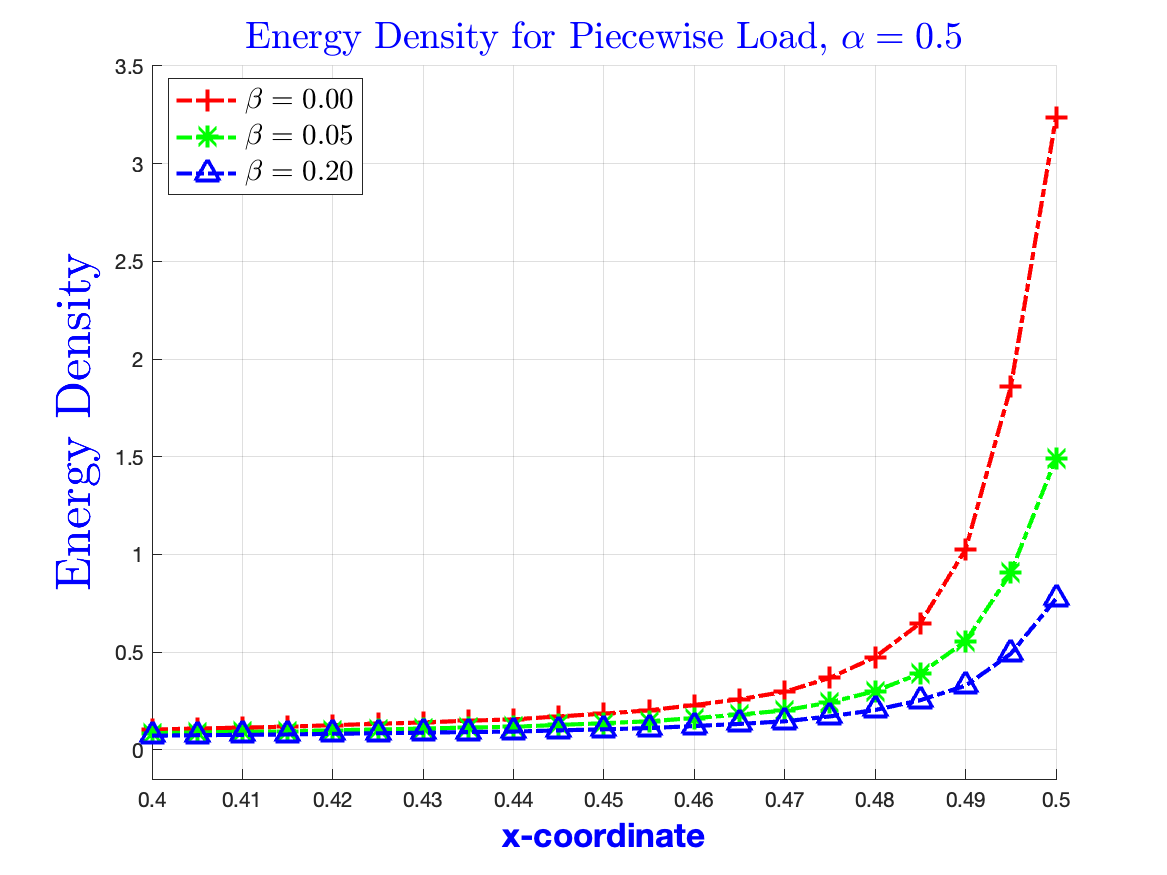}
\caption{Parametric study of the strain-limiting parameter $\beta$ for fibers aligned with the $x$-axis (Case I), with fixed $\alpha=0.5$. From left to right: Variation of Strain, Stress, and Strain Energy Density near the crack tip. Increasing $\beta$ systematically attenuates the magnitude of all crack-tip fields compared to the linear elastic limit ($\beta=0$).\label{fig:mt1}}
\end{figure}

\selectlanguage{american}%
\end{center}

In figure (\ref{fig:mt1}), the model parameter \textbf{$\beta$ }governs
the strain limiting behavior. In the limiting case, $\beta=0,$ the
computed crack-tip fields---strain, stress, and strain energy density---attain
their largest values and essentially reproduce the classical LEFM
response. As $\beta$ increases, the model exhibits a clear attenuation
of all crack-tip quantities, with strain and energy density notably
reduced in the immediate vicinity of the tip. From a physical perspective,
larger value of $\beta$ corresponds to materials that strongly resist
further deformation after a certain point. Under the present piecewise-slope
boundary conditions, this strain-limiting mechanism manifests primarily
as a suppression of the crack opening displacement and a diminished,
highly localized deformation regions surrounding the tip. \begin{figure}[H]
\centering
\includegraphics[scale=0.25]{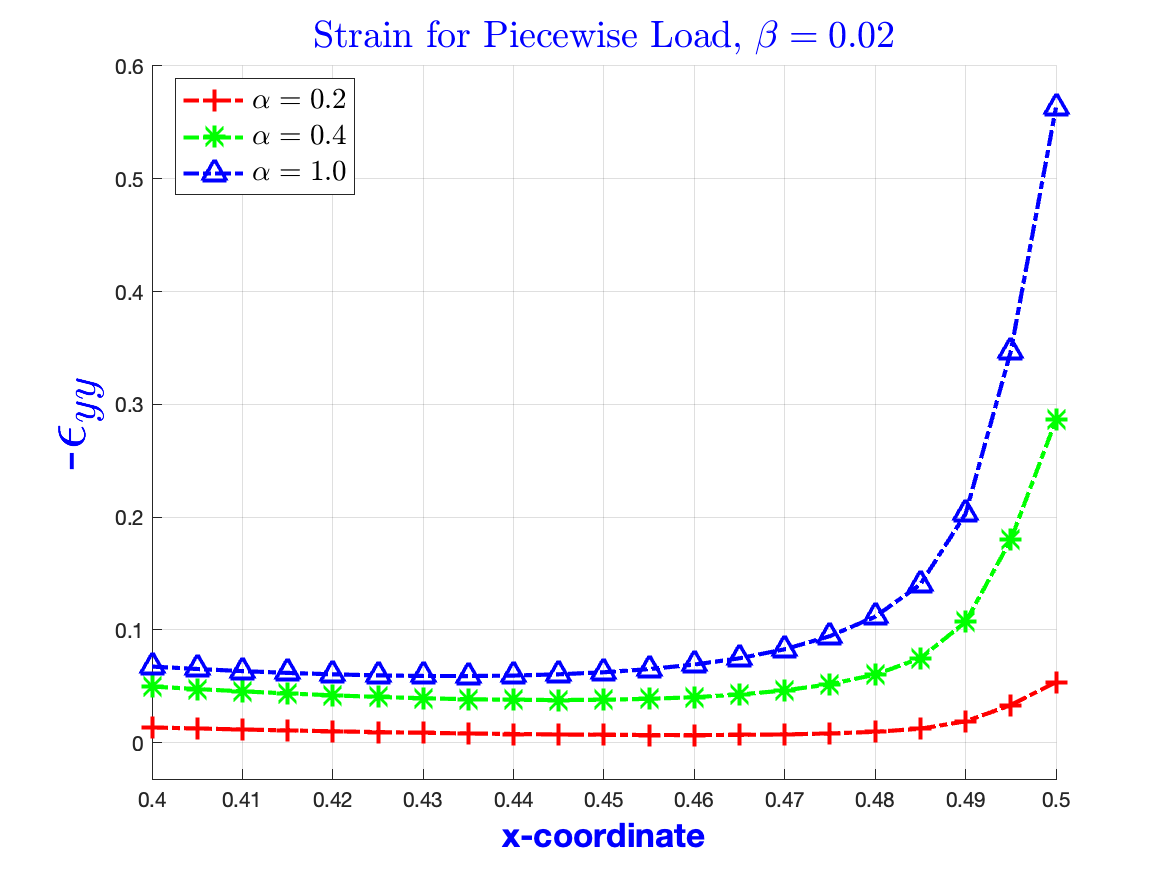}\includegraphics[scale=0.25]{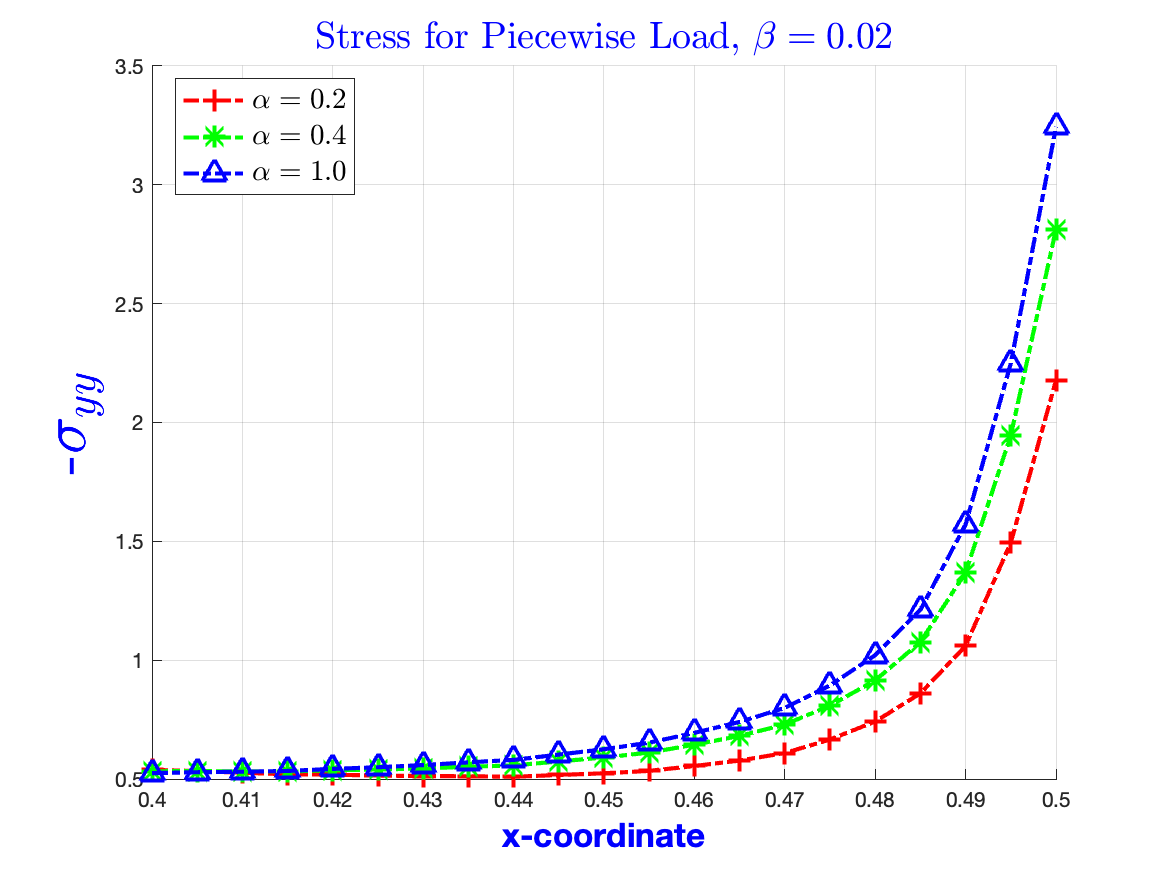}\includegraphics[scale=0.25]{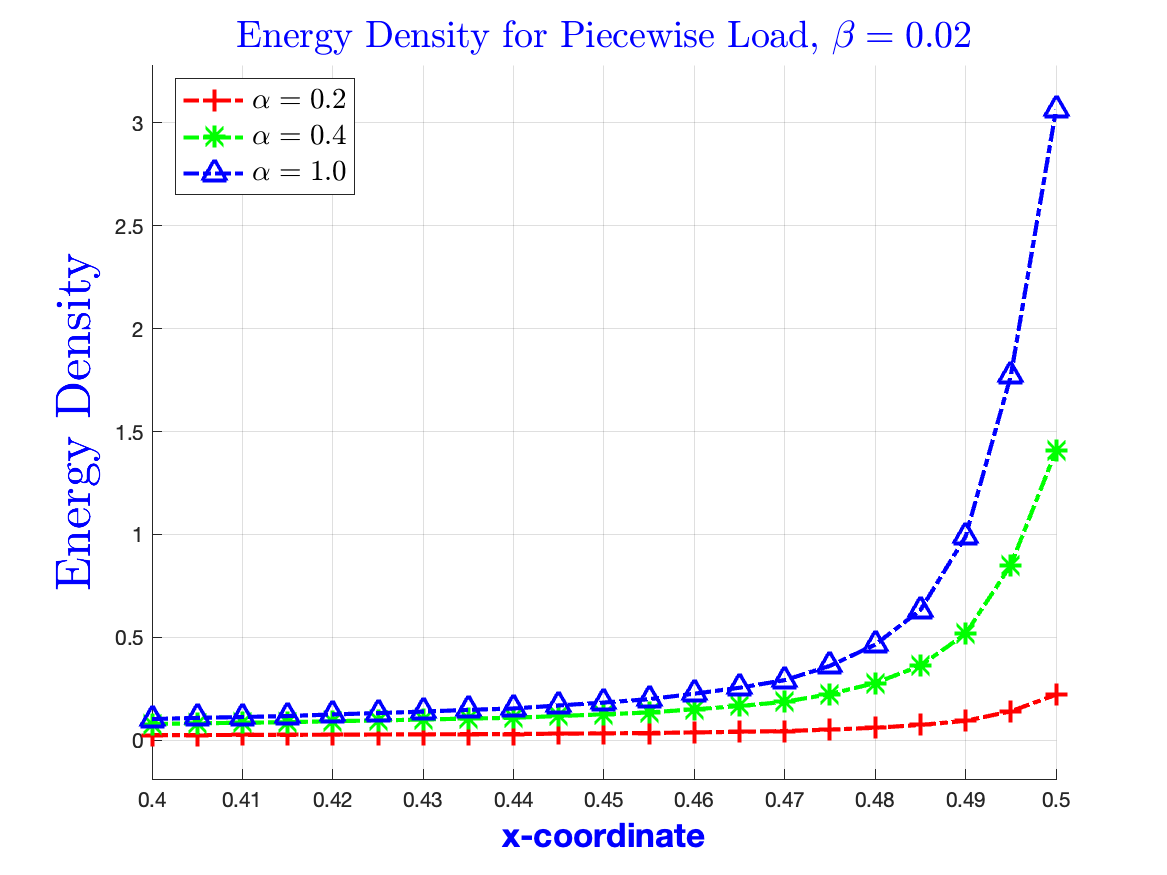}
\caption{Influence of the parameter $\alpha$ on the crack-tip fields for fibers aligned with the $x$-axis (Case I), with fixed $\beta=0.02$. Higher values of $\alpha$ lead to a marked amplification of stress, strain, and strain energy density, indicating sharper gradients and reduced strain-limiting effectiveness near the tip.\label{fig:mt2}}
\protect\selectlanguage{american}%
\end{figure}

\selectlanguage{english}%
Figure (\ref{fig:mt2}) illustrates the influence of the parameter
$\alpha$. In contrast to the behavior observed for $\beta$ , increasing
$\alpha$ leads to a pronounced amplification of the computed fields:
the stress, strain, and strain energy density all grow significantly
in the immediate vicinity of the crack tip. Consequently, materials
characterized by high $\alpha$ are more susceptible to the development
of severe crack-tip fields and are therefore less desirable from a
fracture-resistance perspective, as they exhibit a stronger tendency
toward localized deformation and potential failure. \selectlanguage{american}%
\begin{center}

\begin{figure}[H]
\selectlanguage{english}%
\centering
\includegraphics[scale=0.2]{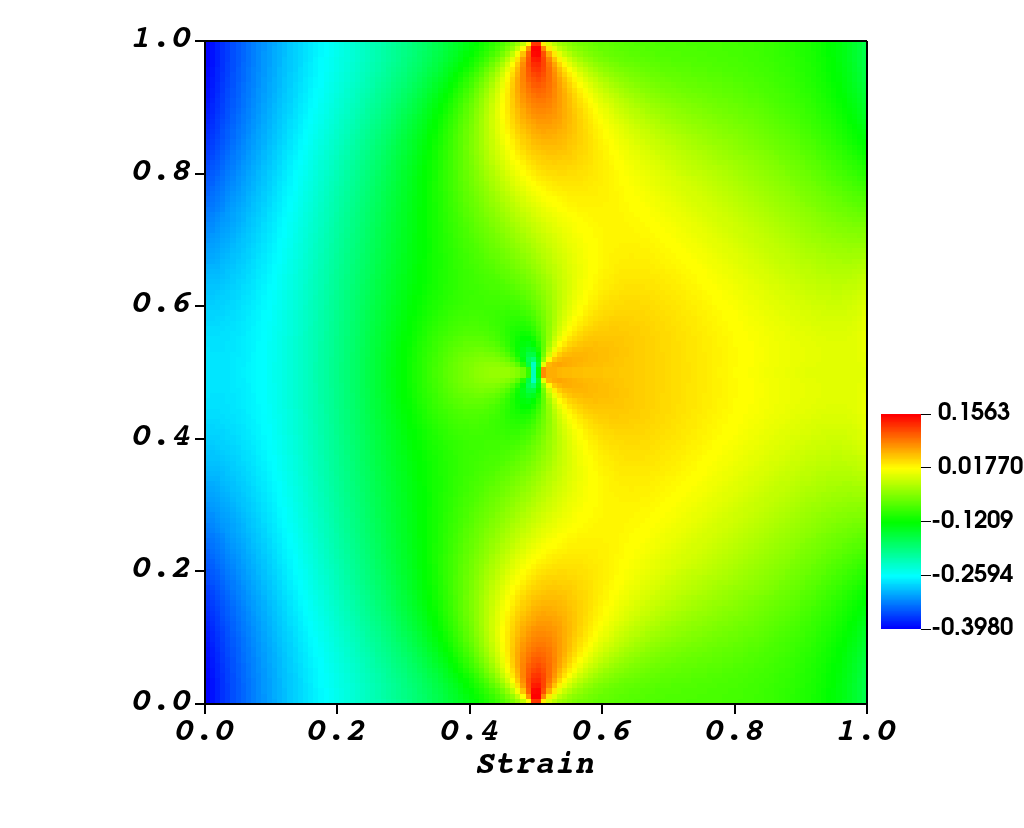}\includegraphics[scale=0.2]{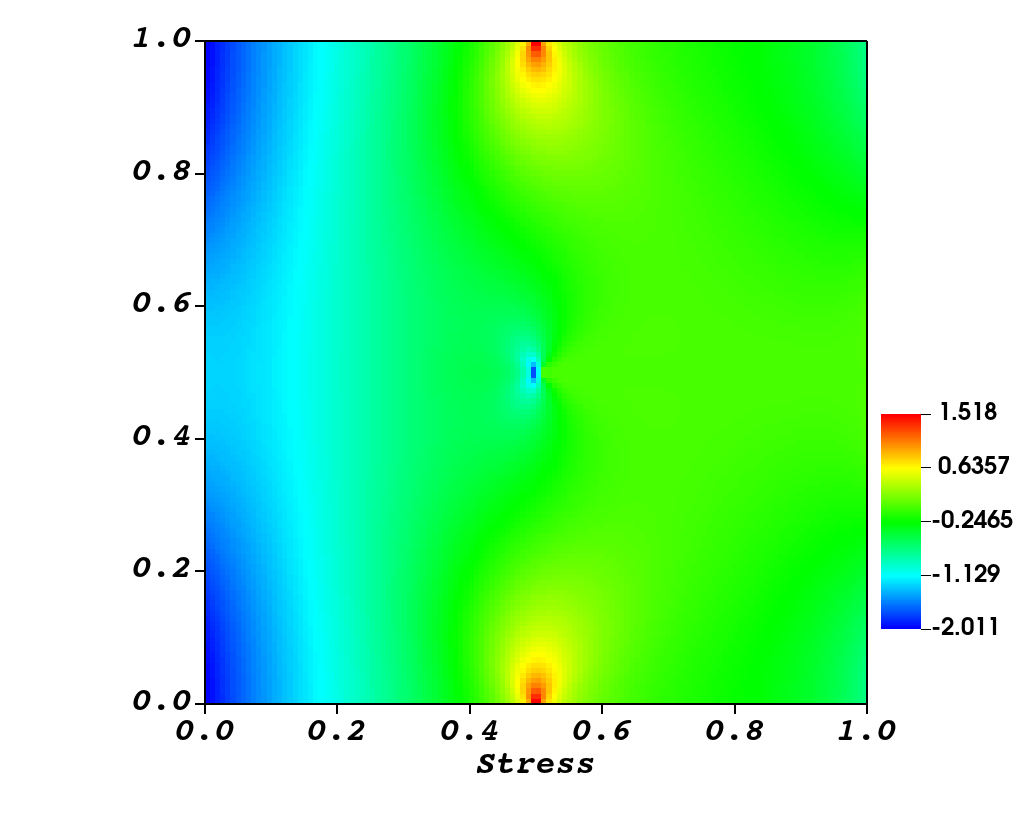}
\caption{Contour plots of the axial strain (left) and axial stress (right) for Case I ($\alpha=0.5, \beta=0.05$). The fields exhibit intensification near the mid-points of the top and bottom boundaries due to the slope discontinuity, along with a localized butterfly pattern surrounding the crack tip.\label{fig:vs1}}
\selectlanguage{american}%
\end{figure}
\end{center}

Figure (\ref{fig:vs1}) illustrates the axial stress and strain of
the domain \foreignlanguage{english}{for $\alpha=0.5$ and $\beta=0.05$.
The most prominent features in both plots are the intensification
of crack tip fields near the mid points of the top and bottom boundaries
of the domain, where the boundary conditions generate the largest
positive magnitudes of both stress and strain (red regions).}

\selectlanguage{english}%
Around the crack tip, the fields are strongly perturbed: a narrow
zone of opposite sign (blue) appears at the crack tip itself, surrounded
by a small halo of moderate values, indicating a local stress--strain
reversal caused by the crack opening. The strain field develops a
characteristic \textquotedblleft butterfly\textquoteright\textquoteright{}
pattern about the crack, whereas the stress field remains more sharply
focused in a small region around the crack tip. Away from these critical zones---particularly in the right half of
the domain---the contours become smoother which indicates a relatively
mild variation of both stress and strain compared to the highly localized
response near the middle of top and bottom edges and the crack tip.
Near the left boundary, both fields show predominantly negative (blue)
values, reflecting the global compressive state induced by the applied
piecewise load. \selectlanguage{american}%

\begin{figure}[H]
\centering
\includegraphics[scale=0.2]{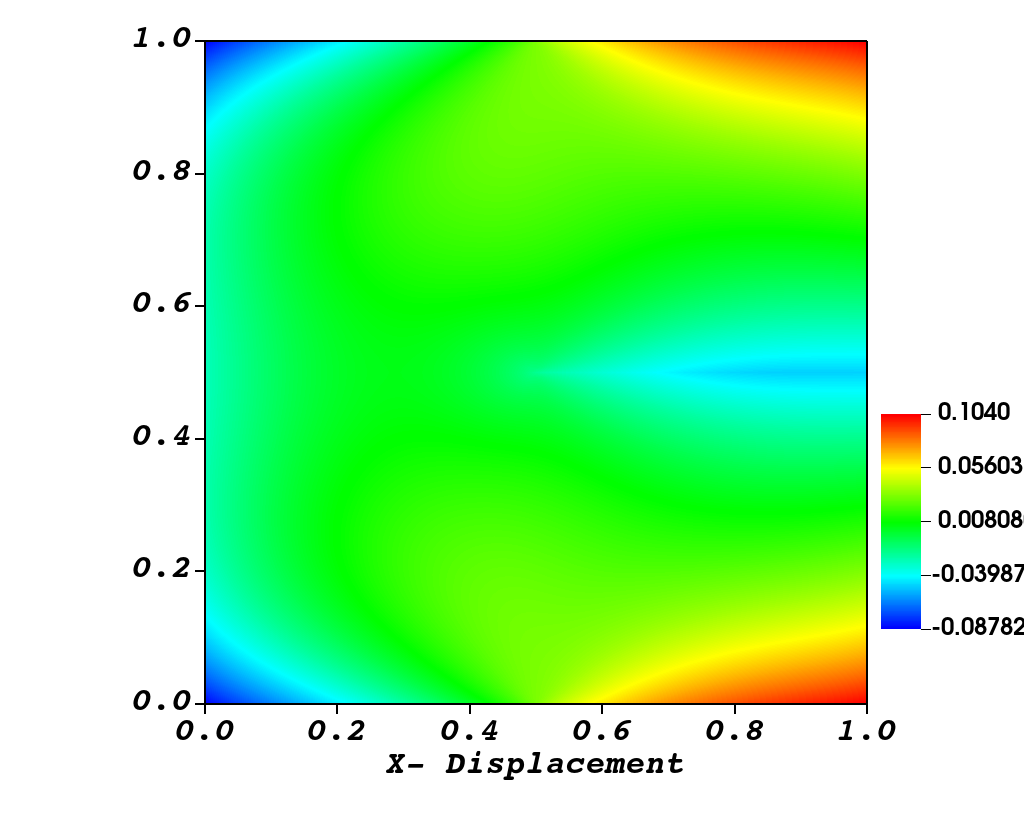}\includegraphics[scale=0.2]{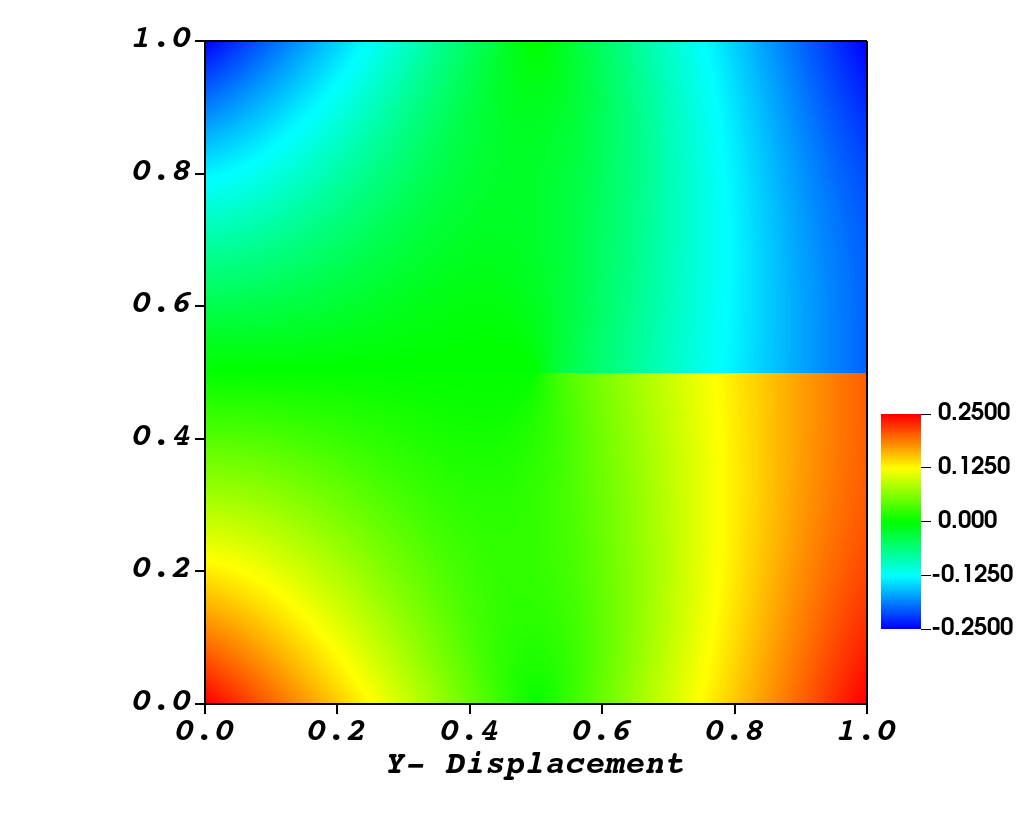}
\protect\caption{Displacement contours for Case I ($\alpha=0.5, \beta=0.05$). Left: Horizontal displacement $u_{1}$, showing mild variation. Right: Vertical displacement $u_{2}$, displaying a sharp antisymmetric transition across the horizontal mid-line induced by the piecewise boundary loading.\label{fig:vs2}}
\protect\selectlanguage{american}%
\end{figure}

Figure (\ref{fig:vs2}) presents the contour plots of the horizontal
displacement $(u_{1})$, and vertical displacement ($u_{2})$. The
horizontal displacement field varies smoothly across the domain, with
the maximum and minimum values occurring near the corners, respectively. A faint disturbance is visible along
the horizontal line passing through the crack, where a narrow bluish
zone indicates a mild inward shift of material adjacent to the crack
faces. The vertical displacement exhibits a much more structured pattern:
the lower half of the domain predominantly undergoes downward motion
(positive displacement), while the upper half undergoes upward (negative
displacement), leading to a sharp transition in displacement values
at mid-height. This reflects the imposed compression and the tendency
of the material to deform symmetrically about the horizontal axis. \begin{figure}[H]
\selectlanguage{english}%
\centering
\includegraphics[scale=0.25]{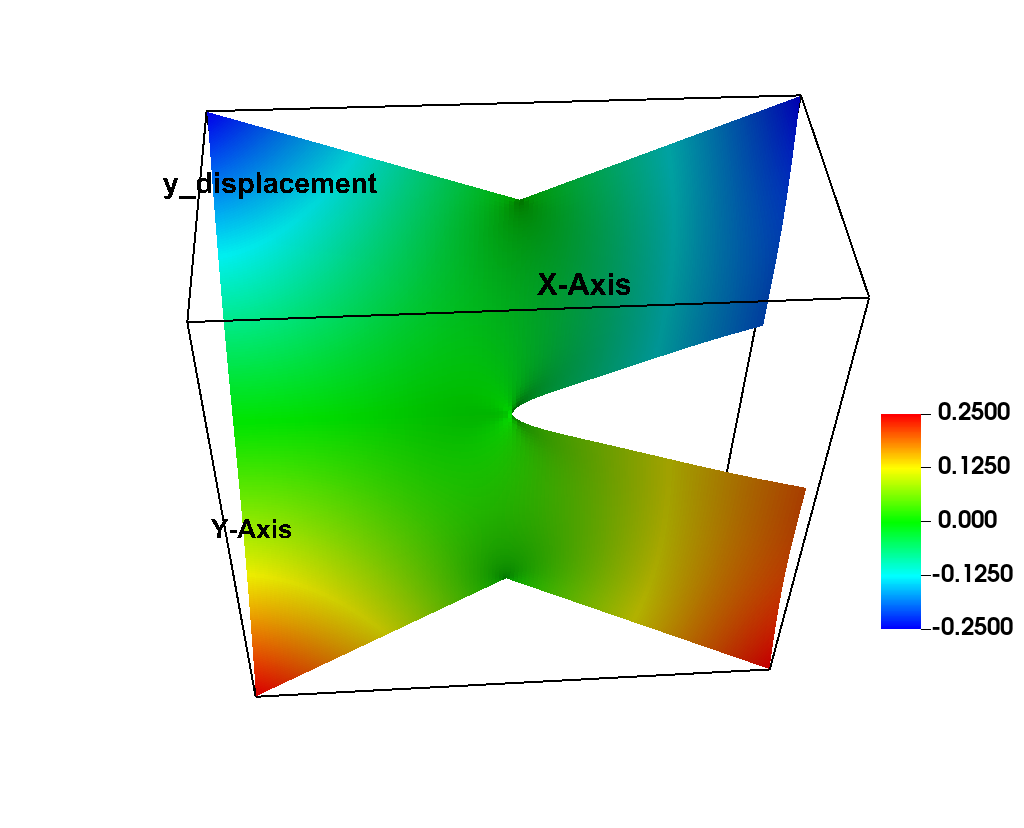}
\caption{Elevated surface view of the vertical displacement $u_{2}$ for Case I. The plot highlights the antisymmetric nature of the deformation, with simultaneous lifting and descending on opposite sides of the crack, and the progressive decay of the bending mode away from the tip.\label{fig:vs3}}
\selectlanguage{american}%
\end{figure}
\foreignlanguage{english}{Figure (\ref{fig:vs3}) shows three-dimensional
elevated view of the vertical displacement $u_{2}$. The antisymmetric
nature of the imposed piecewise-slope boundary condition is clearly
captured: the elevated displacement lifts on one side of the crack
while simultaneously descending on the other, producing a visually
evident sign change in $u_{2}.$ The continuous transition of colors
on the surface reflects a gradual decay in displacement away from
the crack, indicating that although bending remains the dominant mode
of deformation, its influence diminishes progressively toward the
right boundary.}
\selectlanguage{english}%

\subsection{Case II: Fiber Directions along Y axis}

In this case, we assume that the fibers are aligned perpendicular
to the crack. So the structural tensor becomes $M=\boldsymbol{e_{2}}\otimes\boldsymbol{e_{2}}.$
Here $\boldsymbol{e_{2}}=(0,1)$ is unit vector along $y$ - axis,
which implies that the material is stiffer in the direction perpendicular
to the crack. \selectlanguage{american}%
\begin{center}

\selectlanguage{english}%
\begin{figure}[H]
\centering
\includegraphics[scale=0.25]{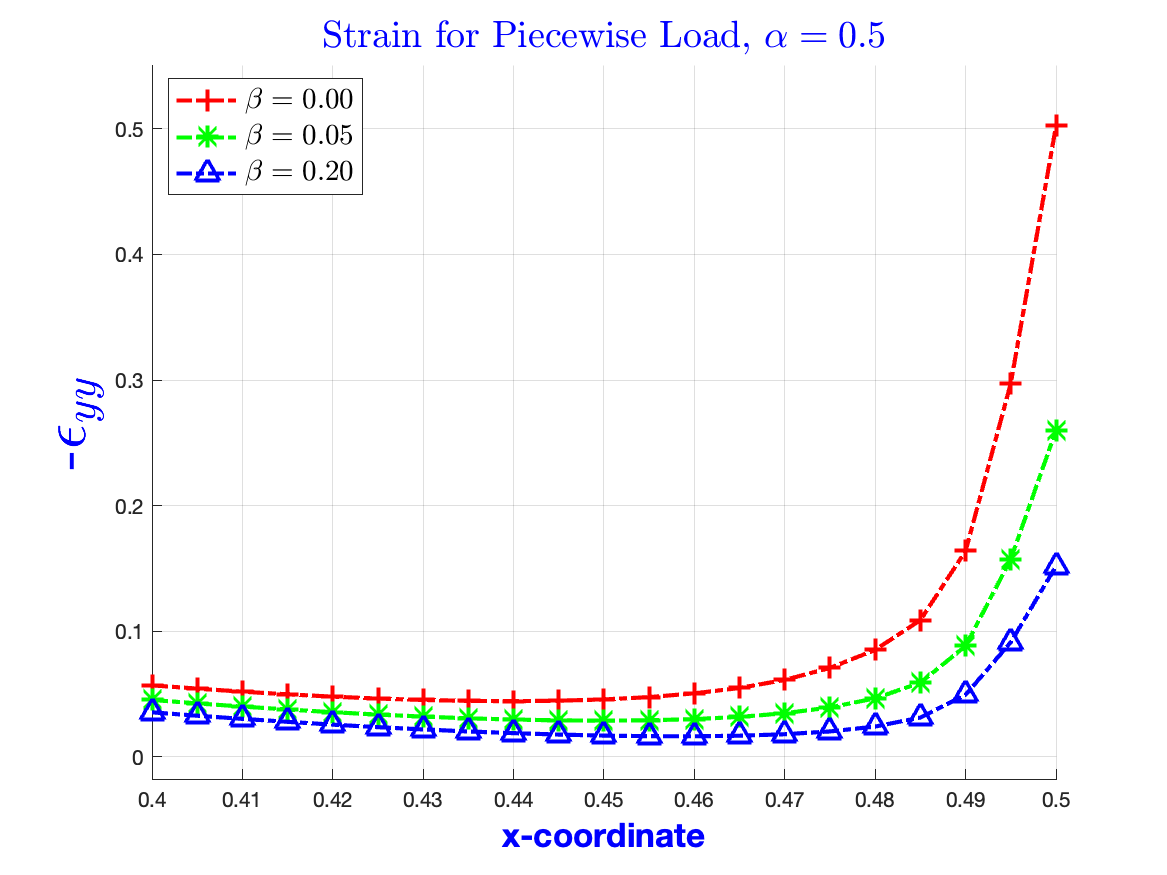}\includegraphics[scale=0.25]{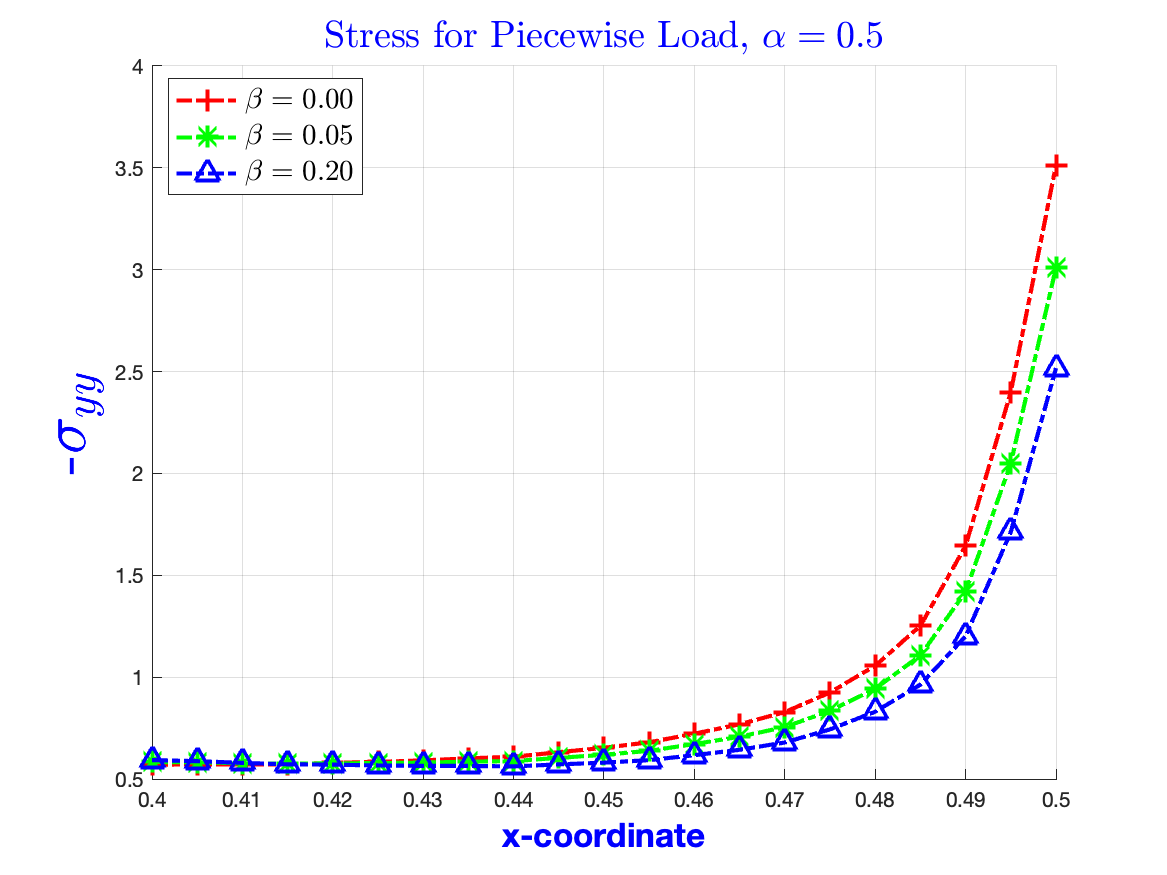}\includegraphics[scale=0.25]{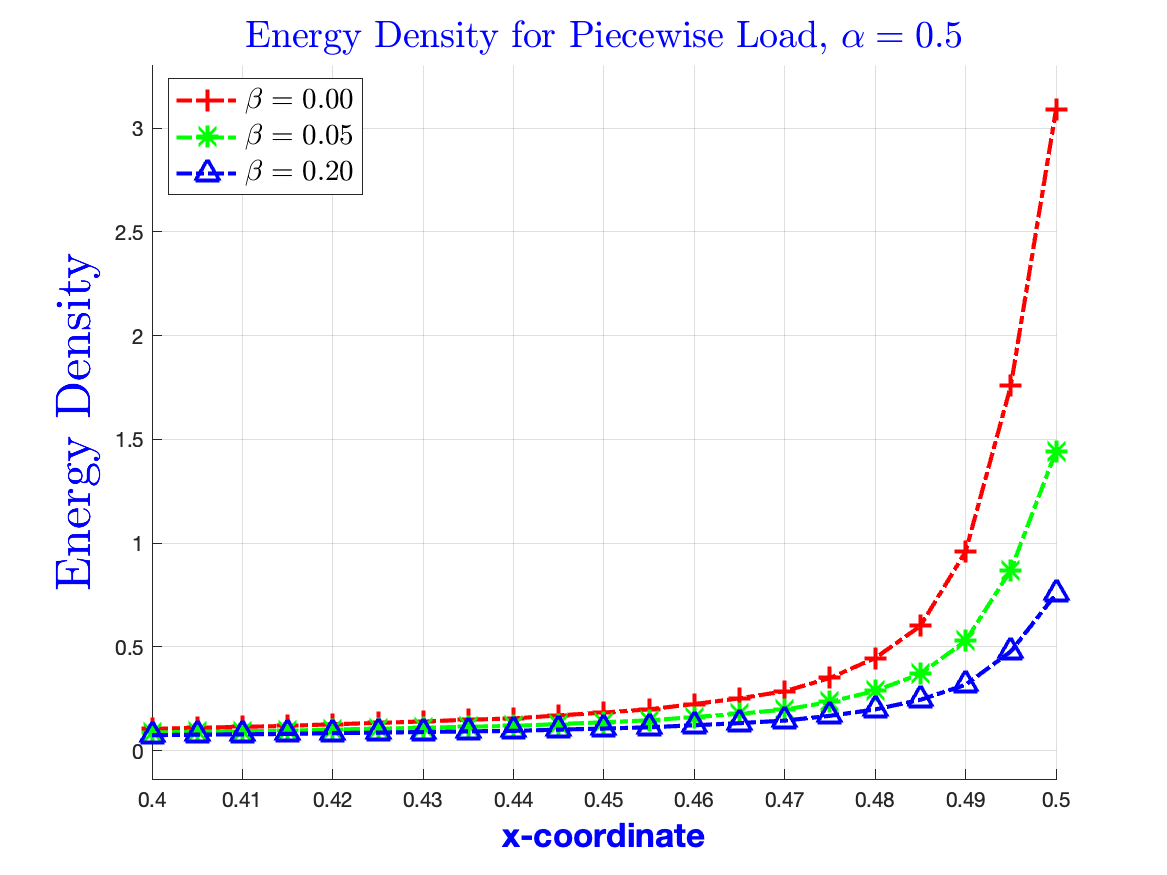}
\caption{Influence of the parameter $\beta$ on the crack-tip fields for fibers aligned with the $y$-axis (Case II). Similar to Case I, increasing $\beta$ significantly suppresses the peak magnitudes of stress, strain, and energy density, reflecting the material's enhanced resistance to crack opening.\label{fig:mt3}}
\end{figure}

\selectlanguage{american}%
\end{center}

\selectlanguage{english}%
Figure~(\ref{fig:mt3}) shows the effect of varying $\beta$ when
the fibers are aligned with the $y$-axis. In this configuration, the
preferred direction is essentially perpendicular to the crack plane,
so the toughening acts directly against the crack opening. As $\beta$
increases, the computed crack-tip fields are systematically suppressed:
the peak values of stress, strain and the strain energy density all
decrease, and the highly stressed zone becomes more confined around
the tip. Mechanically, larger values of $\beta$ in this orientation
strengthen the material\textquoteright s capacity to restrict crack
propagation. The fiber reinforcement aligned with the $y$ axis takes
over a significant share of the imposed load, and curtails additional
crack opening once a certain stress level is reached. \selectlanguage{american}%

\begin{figure}[H]
\centering
\includegraphics[scale=0.25]{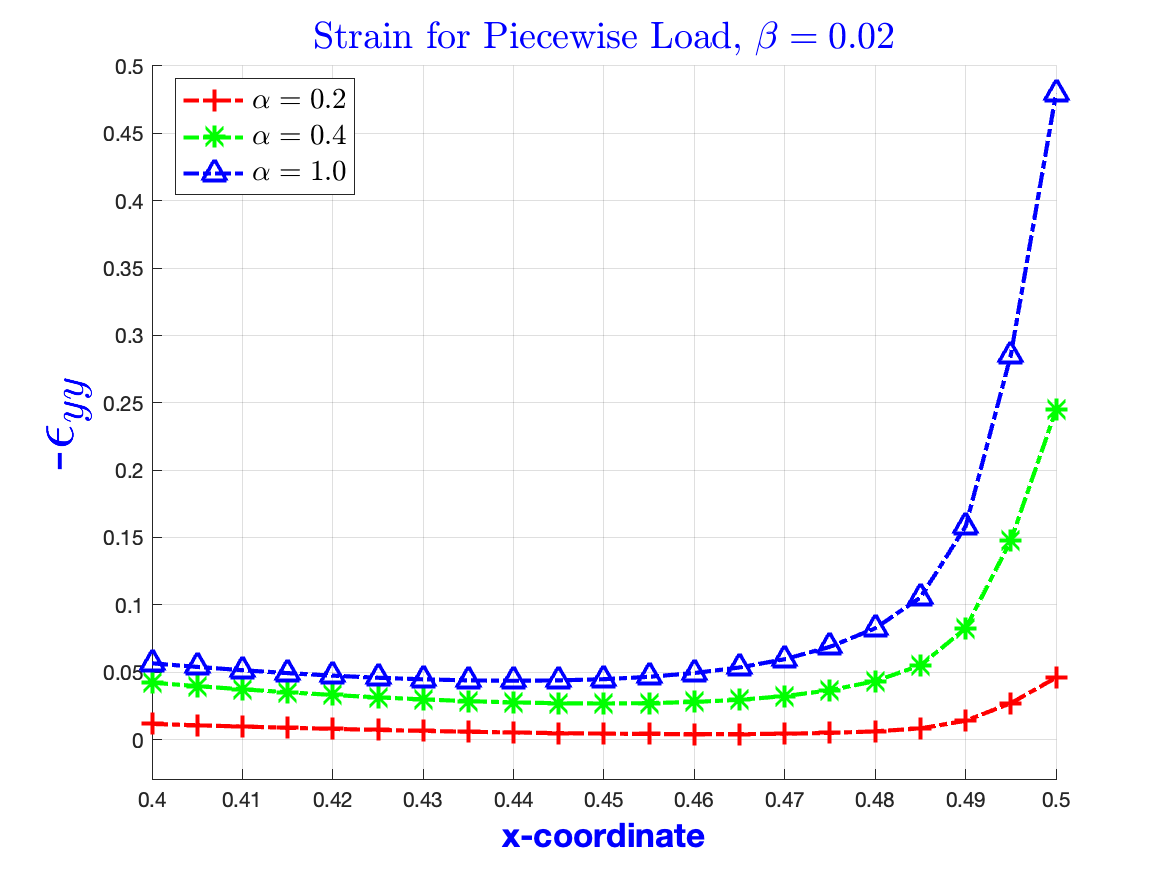}\includegraphics[scale=0.25]{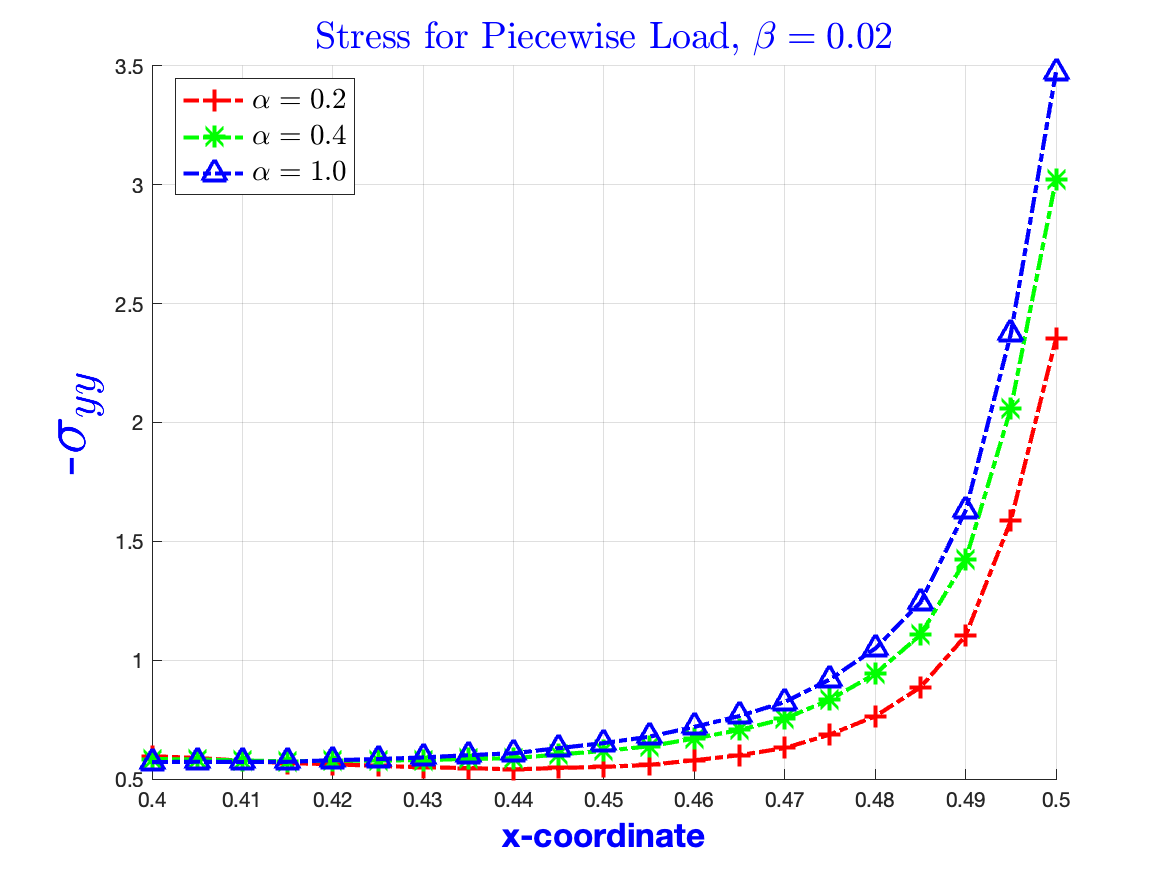}\includegraphics[scale=0.25]{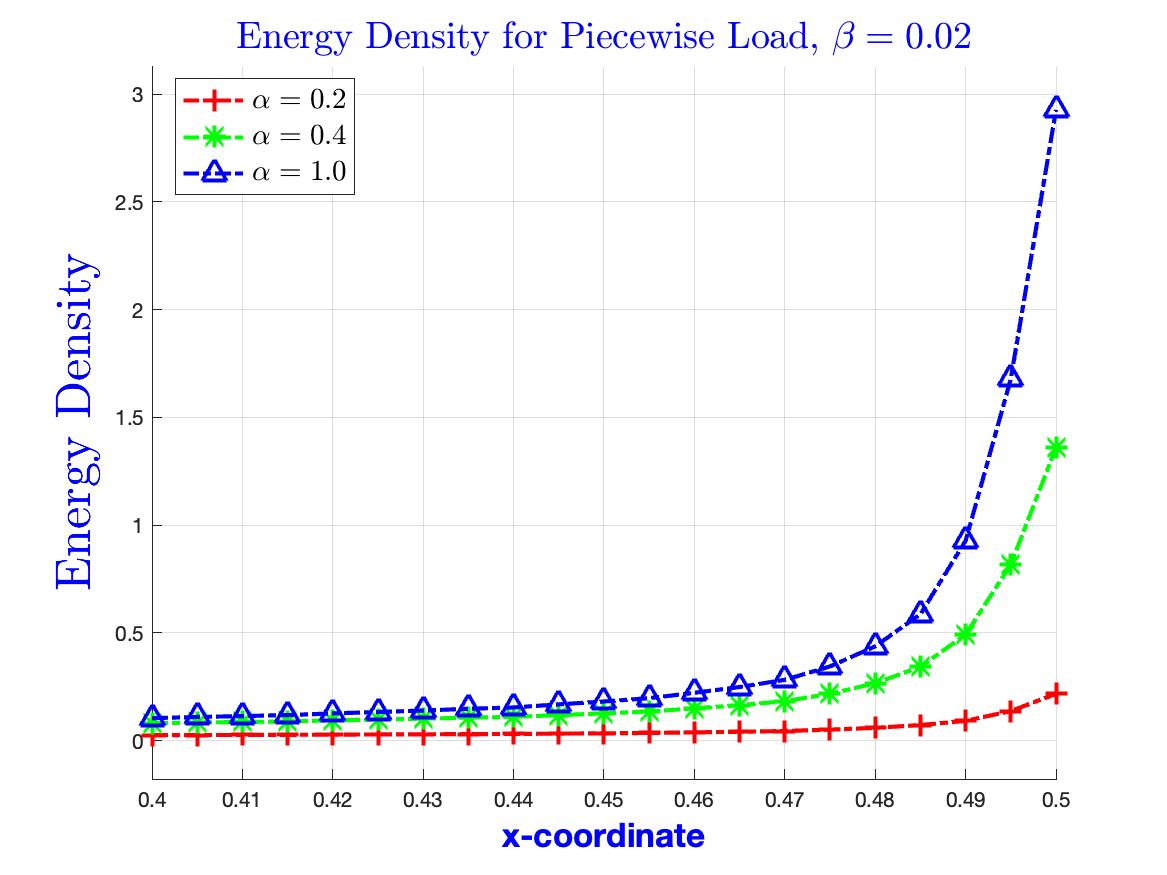}
\protect\caption{Influence of the parameter $\alpha$ on the crack-tip fields for fibers aligned with the $y$-axis (Case II), with $\beta=0.02$. Consistent with Case I, higher $\alpha$ values lead to amplified crack-tip quantities and a reduction in the effective strain-limiting behavior.\label{fig:mt4}}
\protect\selectlanguage{american}%
\end{figure}

\selectlanguage{english}%
Figure~(\ref{fig:mt4}) presents the influence of $\alpha$ for the
case of fibers aligned with the $y$-axis. For fixed $\beta$, higher
values of $\alpha$ lead to a marked amplification of the crack-tip
fields: the stress, strain, and strain energy density all rise in
the vicinity of the tip. This growth in near-tip quantities indicates
a reduction in the effective strain-limiting action of the model,
and consequently a diminished resistance to fracture, since larger
crack-tip stresses and energies generally correlate with an increased
tendency for crack propagation. \selectlanguage{american}%
\begin{figure}[H]
\selectlanguage{english}%
\centering
\includegraphics[scale=0.2]{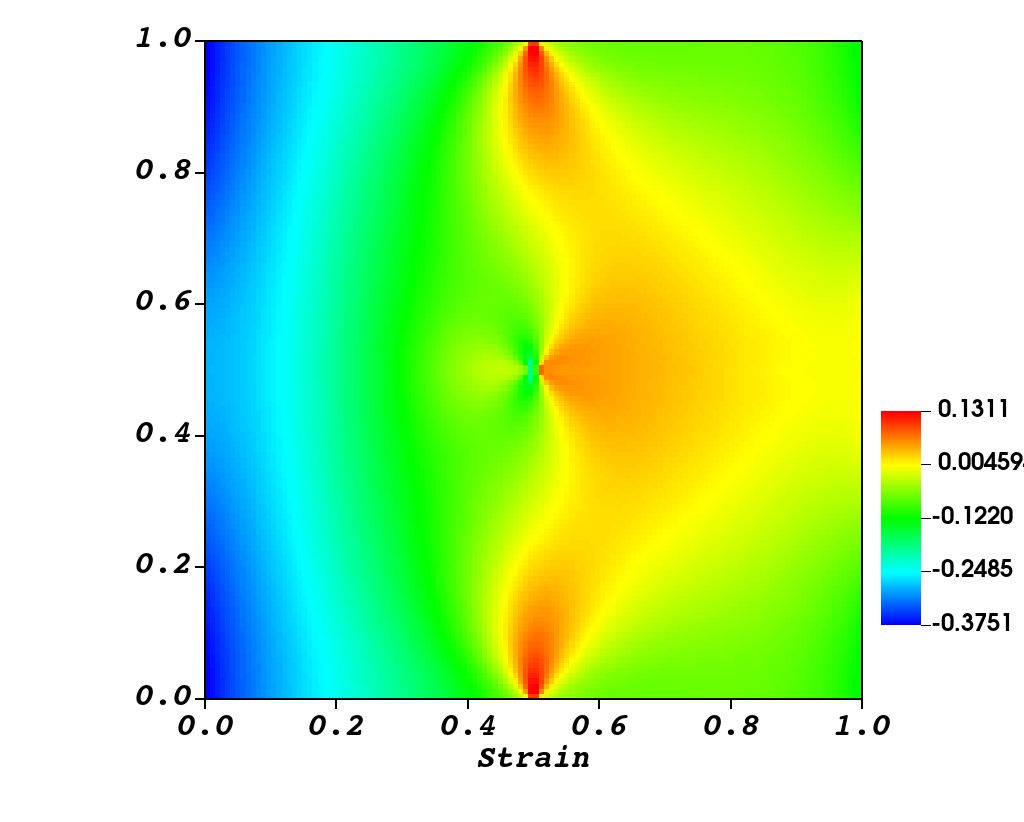}\includegraphics[scale=0.2]{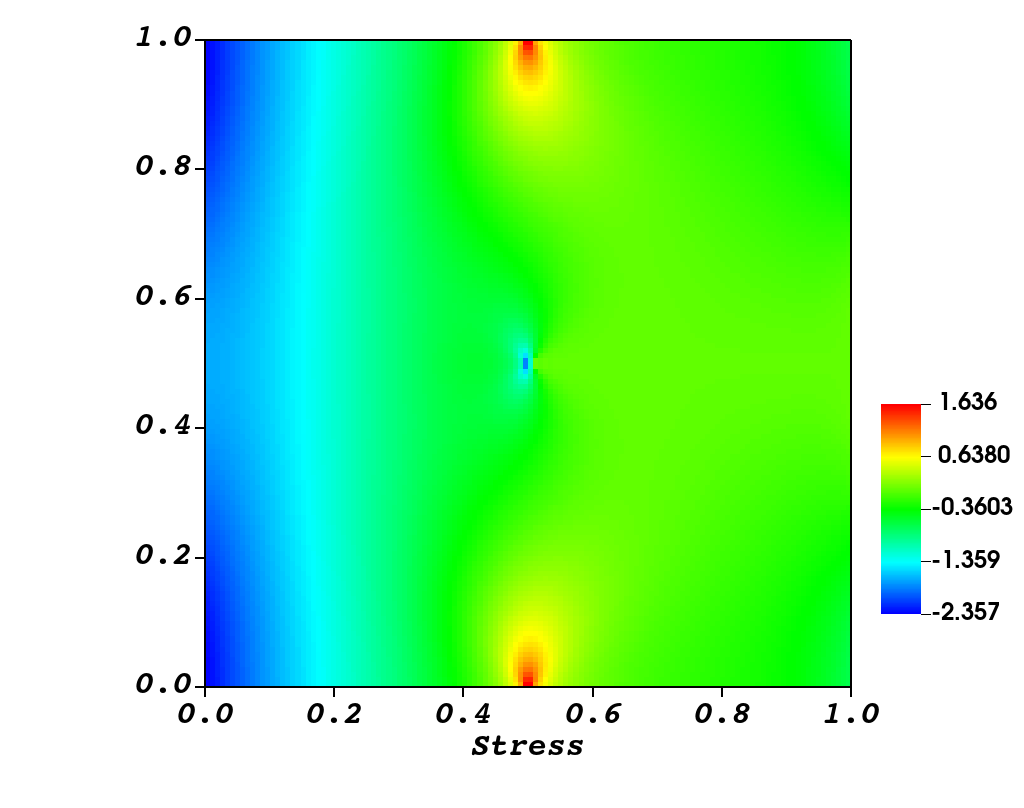}
\caption{Contour plots of axial strain (left) and axial stress (right) for Case II ($\alpha=0.5, \beta=0.05$). The stiffer vertical fiber orientation leads to a vertically elongated strain pattern and strongly localized stress concentrations near the crack tip.\label{fig:vs4}}
\end{figure}

Figure (\ref{fig:vs4}) shows the axial stress and strain distributions
in the domain for $\alpha=0.5$ and $\beta=0.05$ when the fibers
are oriented along the $y$-axis. In this case, the most striking
features appear along the upper and lower boundaries, where the imposed
displacement produces pronounced vertical stretching and compression.
These regions exhibit the highest positive intensity in both fields
(red zones), reflecting the directional stiffness induced by the fiber
orientation. Near the crack tip, both fields change abruptly. A small blue region
appears exactly at the tip, showing that the stress and strain switch
sign there. This region is surrounded by a smoother band of moderate
values which indicates a localized reversal and redistribution of
stress and strain generated by the crack opening. The strain contours
extend more noticeably in the vertical direction, creating a vertically
elongated pattern, while the stress field remains more concentrated
and strongly localized around the crack tip due to the anisotropic
response. Farther from the crack and the loading boundaries---particularly
across the right portion of the domain---the contours gradually flatten
and become more uniformly spaced, signifying weaker gradients and
a relatively undisturbed deformation state. Along the left boundary,
both fields take predominantly negative (blue) values, consistent
with the compressive nature of the applied loading. \begin{figure}[H]
\centering
\includegraphics[scale=0.2]{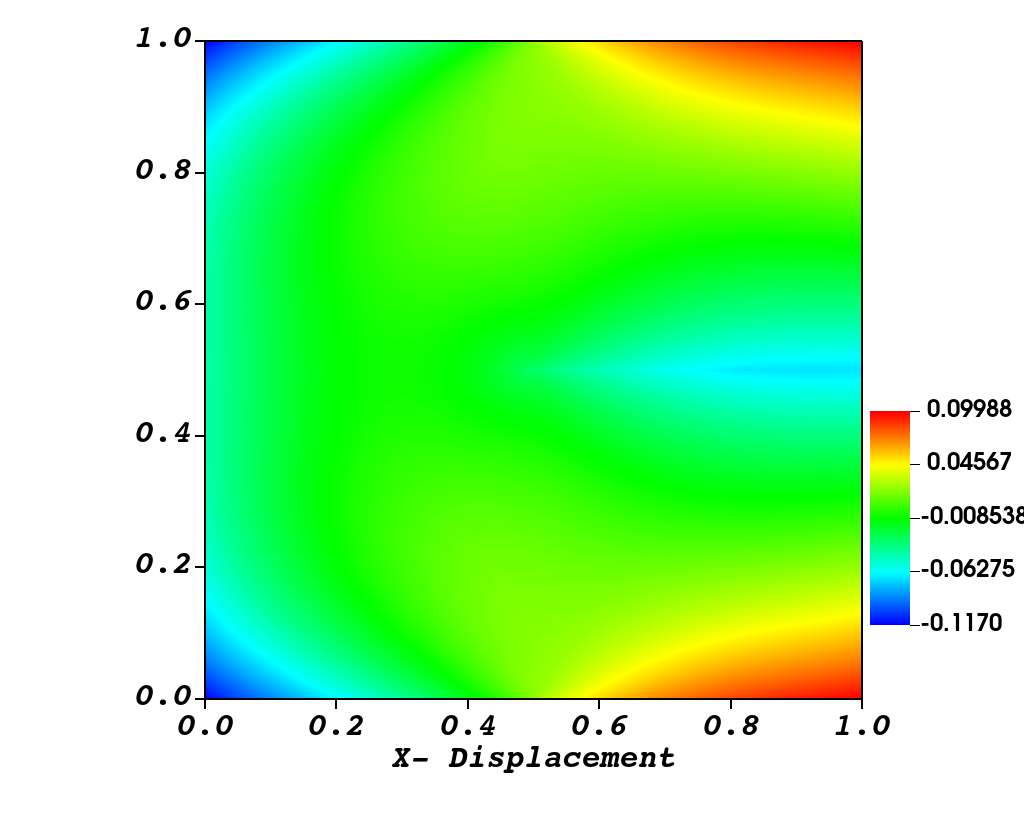}\includegraphics[scale=0.2]{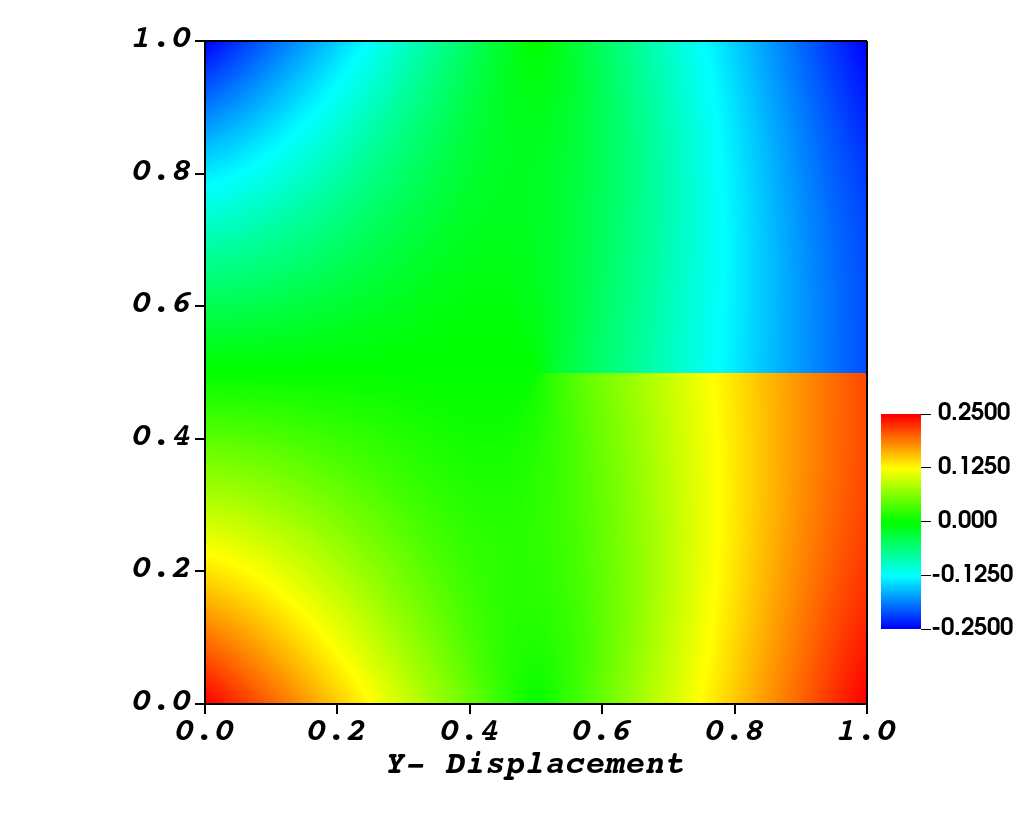}
\protect\caption{Displacement contours for Case II ($\alpha=0.5, \beta=0.05$). Left: Horizontal displacement $u_{1}$. Right: Vertical displacement $u_{2}$. The vertical displacement field is dominant and exhibits a clear sign reversal across the mid-height, driven by the boundary loading.\label{fig:vs5}}
\protect\selectlanguage{american}%
\end{figure}

Figure (\ref{fig:vs5}) displays the contour plots of the horizontal
displacement $\left(u_{1}\right)$ and the vertical displacement $(u_{2})$
for the case in which the fibers are aligned with the $y$-axis. The
horizontal displacement $u_{1}$ is considerably smaller in magnitude
and varies smoothly throughout the plate. A gentle diagonal gradient
is observed, with negative values near the top-left corner and positive
values toward the bottom-right corner. A slight disturbance appears
along the horizontal line containing the crack, but its effect on
the global displacement field remains mild. The vertical displacement
field exhibits the dominant deformation pattern with a clear sign
change across the mid-height of the domain. \begin{figure}[H]
\selectlanguage{english}%
\centering
\includegraphics[scale=0.25]{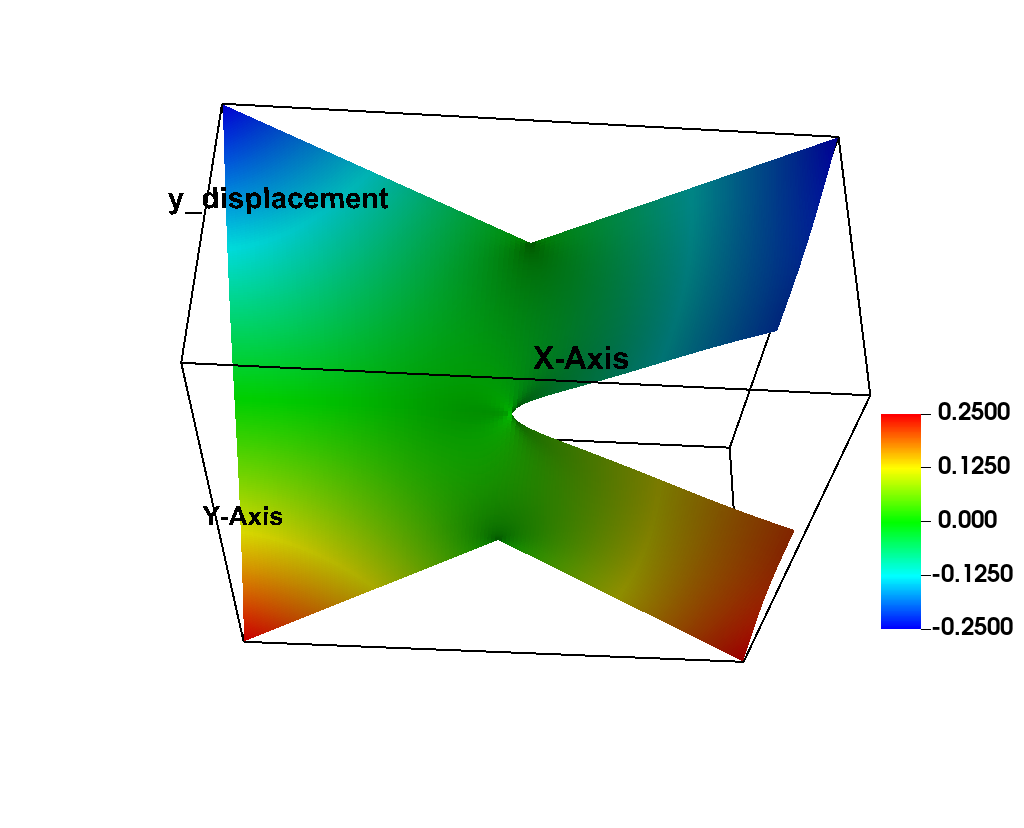}
\caption{Elevated surface view of the vertical displacement $u_{2}$ for Case II. The deformation exhibits a distinctive saddle-like profile with a localized peak near the crack tip, emphasizing the antisymmetric bending response.\label{fig:vs6}}
\selectlanguage{american}%
\end{figure}
The\foreignlanguage{english}{ elevated view of the vertical displacement
$u_{2}$ is presented in Figure (\ref{fig:vs6}). The deformation
profile exhibits a distinctive saddle-like shape, characterized by
a localized elevation near the crack tip and a corresponding downward
deflection along the opposite half of the plate. The peak vertical
displacement occurs at the crack tip, producing a clear antisymmetric
bending response across the mid-span. The transition from negative
to positive displacement is smooth, captured by the continuous color
gradient, and reflects the strong coupling between the loading configuration
and the material anisotropy.}
\selectlanguage{english}

\section{Conclusion}

We have studied crack-tip fields in a transversely isotropic strain-limiting
elastic body subjected to a class of non-uniform boundary conditions
in the form of piecewise linear slope-type displacements. The mechanical
response was modeled within strain-limiting elasticity, in which the
linearized strain is related to the Cauchy stress through a nonlinear,
algebraic constitutive law that enforces a bounded strain response.
The resulting quasi-linear elliptic boundary value problem was obtained
in weak form and approximated by a conforming finite element method
combined with a Picard linearization scheme. The numerical experiments reveal how the model parameters fundamentally
affect crack-tip fields in transversely isotropic materials. The parameter
$\beta$ acts as the principal strain-limiting factor: increasing
$\beta$ systematically diminishes the magnitudes of the crack-tip
stress, strain, and strain-energy density, effectively moderating
the deformation in the vicinity of the crack. This behavior accords
with the constitutive structure of the model---larger $\beta$ enforces
a more strict bound on admissible strain levels---thereby reducing
the mechanical driving force available for crack propagation. In this
regard, $\beta$ may be interpreted as introducing a material-level
toughening mechanism that delays the onset of catastrophic fracture.
This damping influence persists irrespective of fiber orientation,
although its manifestation becomes more pronounced when the preferred
direction is orthogonal to the crack, since the fibers then act directly
against displacement. Conversely, the parameter $\alpha$ plays an opposite role. The simulations
show that increasing $\alpha$ leads to sharper amplification of the
crack-tip fields---stress, strain, and strain-energy density all
rise significantly as $\alpha$ increases. Materials characterized
by large $\alpha$ respond more aggressively to imposed loads, producing
steeper gradients near the crack and thereby exhibiting lower resistance
to fracture initiation. In both fiber orientations, larger $\alpha$
values make the deformation much more concentrated near the crack
tip. This leads to sharper gradients in stress and strain and an increase
in stored strain energy around the crack. Such strong localization
is characteristic of brittle materials, where cracks tend to propagate
suddenly once the energy concentration exceeds a critical threshold.
These findings demonstrate that the interplay between $\alpha$, $\beta$, and fiber orientation offers a powerful parametric
mechanism for tailoring crack-tip behavior in strain-limiting materials.
By appropriately selecting $\beta$, one may suppress
non-physical strain growth, while careful control of $\alpha$
enables tuning of the model's sensitivity to stress
escalation near defects. The results of this study therefore extend
the existing body of work on strain-limiting elasticity by showing
that non-uniform, piecewise slope loading not only provides a realistic
boundary condition class, but also exposes distinct fracture characteristics
that are invisible under uniform loading. This paves the way for future
investigations involving mixed-mode fracture and phase-field approaches
under non-smooth load distributions.

\section*{Acknowledgement}
SMM thanks the  National Science Foundation for its financial support under Grant No. 2316905. SG would like to thank the University of Texas Rio Grande Valley for providing a Presidential Research Fellowship during his PhD studies. 

\bibliographystyle{plain}  
\bibliography{p5_references.bib}

\end{document}